\documentclass[12pt,a4paper]{article}
\usepackage{amsfonts,amsmath,amssymb,amscd}
\usepackage{fancybox}
\usepackage[latin1]{inputenc}
\usepackage[usenames]{color}
\usepackage{graphicx}
\usepackage{float}
\usepackage{fullpage}
\usepackage{epsfig,amssymb}
\usepackage{enumerate}
\usepackage{subcaption}

\definecolor{avocado}{rgb}{0.34,0.51,0.01}

\newcommand{\intprod}{\mathbin{\raisebox{\depth}{\scalebox{1}[-1]{$\lnot$}}}}

\newtheorem{thm}{Theorem}[section]

\newtheorem{lem}[thm]{Lemma}
\newtheorem{prop}[thm]{Proposition}

\newtheorem{rem}[thm]{Remark}

\begin{document}
	
\title{\textbf{On p-Willmore Disks with Boundary Energies}}

\author{Anthony Gruber, \'Alvaro P\'ampano and Magdalena Toda}

\date{}

\maketitle 

\begin{abstract}
\noindent We consider an energy functional on surface immersions which includes contributions from both boundary and interior.  Inspired by physical examples, the boundary is modeled as the center line of a generalized Kirchhoff elastic rod, while the interior term is arbitrarily dependent on the mean curvature and linearly dependent on the Gaussian curvature. We study  equilibrium configurations for this energy in general among topological disks, as well as specifically for the class of examples known as p-Willmore energies.
\\

\end{abstract}

\section{Introduction}

The significance of variational problems throughout mathematics and the physical sciences cannot be overstated. From surfactant films and lipid membranes to material interfaces (to name just a few), equilibrium solutions are invaluable  for the theoretical modeling of  physically-observable quantities.  Consequently, variational equilibria have enjoyed centuries of interest, not only in mathematics but also in fields ranging from biology to architecture.

Originating from the elastic theory pioneered by J. Bernoulli, D. Bernoulli and L. Euler in the case of curves and S. Germain and S. Poisson in the case of surfaces, the study of variational problems belongs mostly to the fields of geometric PDE and the calculus of variations.  Indeed, the bending and twisting of matter in space induces elastic deformations typically modeled through so-called Lagrangians, which depend on the mathematics of curvature. For instance, in 1811 S. Germain suggested using the integral of a symmetric polynomial of even-degree in the principal curvatures to study the physical system associated with an elastic plate.  In particular, the early development of geometric PDE was motivated by similar problems which 
attracted the attention of many outstanding mathematicians including Blaschke, Kirchhoff, Thomsen, Willmore, and others.

Arguably the most studied energy in this mathematical context is the total squared curvature or bending energy. For a surface immersion $X:\Sigma\rightarrow\mathbb{R}^3$, this energy is classically presented as
\begin{equation}\label{nonconfwillmore}
\mathcal{W}[X]:=\int_\Sigma H^2\,d\Sigma\,,
\end{equation}
and represents one of the simplest functionals suggested by S. Germain, see e.g. \cite{Germain}. Later reintroduced by T. Willmore \cite{Willmore}, the bending energy $\mathcal{W}[X]$, \eqref{nonconfwillmore}, is usually referred to as the Willmore energy.

The Willmore energy is conformally invariant among closed surfaces, an essential property shown by W. Blaschke \cite{Blaschke} in the early 20th century. Conversely, for surface immersions with boundary, this conformal invariance is only preserved by considering the so-called conformal Willmore energy
\begin{equation}\label{confwillmore}
\mathcal{W}[X]:=\int_\Sigma\left(H^2-K\right)d\Sigma\,.
\end{equation}
For closed surfaces, the conformal Willmore energy is variationally equivalent to its namesake, since the Gauss-Bonnet Theorem implies that the total Gaussian curvature is a topological invariant.   Consequently, the same notation has been used here to denote both functionals.

Beyond their mathematical interest, Willmore-type energies have been widely studied in biology as models for naturally occurring objects. In the seminal work of P. Canham \cite{Canham}, the minimization of the energy $\mathcal{W}[X]$, \eqref{confwillmore}, was proposed as a possible explanation for the biconcave shape of red blood cells. Moreover, the main constituents of cellular membranes in most living organisms are lipid bilayers formed from a double layer of phospholipids, which are frequently modeled as mathematical surfaces  since they are very thin compared to the size of the observed cells or vesicles. An early pioneer of such modeling was W. Helfrich \cite{Helfrich}, who proposed an extension of the energy functional $\mathcal{W}[X]$, \eqref{confwillmore}, based on liquid crystallography.  The Helfrich energy has the general form
\begin{equation}\label{Helfrich}
\mathcal{H}[X]:=\int_\Sigma \left(\sigma\left[H-c_o\right]^2 + \eta K\right)d\Sigma\,,
\end{equation}
where $H$ and $K$ are the mean resp. Gaussian curvatures of the surface, while $\sigma>0$ and $\eta$ represent bending rigidities which depend on the composition of the membrane itself. Here the parameter  $c_o$ represents the spontaneous curvature, and its sign favors the local geometry to be spherical, planar or hyperbolic \cite{RL}. Under the simplifying assumption that the membrane is homogeneous, each of $\sigma>0$, $\eta$, and  $c_o$ is a physical constant. Note that any spontaneous curvature $c_o$ originates primarily from the asymmetry between the two layers of the membrane \cite{Seifert}, although it may also arise from differences in the chemical properties of the fluids on both sides of the lipid bilayer \cite{Betal}. Geometrically, this can be measured as the difference in area between the two layers, giving rise to the previous model $\mathcal{H}[X]$, \eqref{Helfrich}, for the Helfrich energy.

In the case of closed surfaces, equilibria and minimizers of the Helfrich energy $\mathcal{H}[X]$, \eqref{Helfrich}, are usually considered under some constraints  motivated by the involved physics. Indeed, phospholipid molecules are composed of a hydrophilic head and a hydrophobic tail. Therefore, when suspended in an aqueous solution, they spontaneously aggregate in order to shield the tails from the solvent, preventing an energetically unfavorable scenario and resulting in a closed bilayer with heads pointing outwards \cite{AB}. Consequently, the total area or total enclosed volume is often fixed during the mathematical analysis of these phenomena. In this setting, the problem of existence of minimizers for the Helfrich energy $\mathcal{H}[X]$, \eqref{Helfrich}, has been approached from different viewpoints. The existence of a minimizer among a suitable class of axially symmetric (possibly singular) surfaces under fixed surface area and enclosed volume constraints was proven in \cite{Choksi} using the direct method of calculus of variations, i.e. the geometric measure theory. More recently in \cite{Mondino}, a parametric approach based on PDE theory and functional analysis was used to prove the existence and regularity of minimizers in the class of weak (possibly branched and bubbled) immersions of the sphere.

On the other hand, the formation of pores in biological membranes is quite common due to different types of stimuli such as mechanical stresses and thermal instabilities, and stabilized open bilipid membranes may exist. This provides the potential for equilibrium membranes with edges \cite{BMF}, which have been studied in \cite{CGS}, \cite{Tu}, \cite{Tu-Ou-Yang} and \cite{Z} (to mention some) through a combination of the Helfrich energy $\mathcal{H}[X]$, \eqref{Helfrich}, and an inelastic boundary line tension.  In fact, it turns out that the line tension of the boundary represents only a simplified model of the real behavior, as the computations in \cite{AB} show the boundary energy includes a general elastic term measuring not only bending but also twisting. Different studies of the Helfrich energy $\mathcal{H}[X]$, \eqref{Helfrich}, with elastic boundary can be found in e.g. \cite{BMF}, \cite{Gnew}, \cite{Palmer-Pampano-2} and the references therein.

Presently, we consider immersions of a compact surface of disk type, $X:\Sigma\rightarrow\mathbb{R}^3$, along with the total energy
\begin{equation}\label{E}
E[X]:=\int_\Sigma P(H)\,d\Sigma+\eta\int_\Sigma K\,d\Sigma+\oint_{\partial\Sigma}\left(\Lambda(\kappa)+\varpi\tau+\beta\right)ds\,,
\end{equation}
for suitable constant energy parameters $\eta$, $\varpi$ and $\beta$ along with analytic functions $P(H)$ and $\Lambda(\kappa)$ which are motivated by the physics. The present theoretical contribution is a mathematical analysis of the energy $E[X]$, \eqref{E}, with special attention given to the differential geometry of equilibria.

By considering an arbitrary surface energy density depending on the mean curvature, we obtain an extension of the first term in the Helfrich energy $\mathcal{H}[X]$, \eqref{Helfrich}, which includes (but is not restricted to) the original suggestions of S. Germain \cite{Germain} as well as the Willmore and p-Willmore energies (see e.g. \cite{Gruber-Toda-Tran-1}). The additional inclusion of a total Gaussian curvature term is  motivated mathematically by the conformal Willmore energy $\mathcal{W}[X]$, \eqref{confwillmore}, and the Helfrich energy $\mathcal{H}[X]$, \eqref{Helfrich}, but also plays an important role from a physical viewpoint.  Indeed, current literature suggests this term as one of the major contributions to the free energy of fusion stalks in phospholipids, as well as the tendency for membranes at phase boundary to form intermediate structures (see e.g. \cite{Siegel-Kozlov}). 

In the same fashion, the present boundary energy extends those proposed before, arising from consideration of the boundary as the center line of a generalized Kirchhoff elastic rod. Determining the shape of ideal elastic rods is a classical problem in mechanics and the calculus of variations (see e.g. \cite{Langer-Singer} and references therein), and its associated energy represents a meaningful generalization of previously suggested boundary energies. In particular, the present energy accounts for the possible twisting of the boundary as suggested in \cite{AB}, while also allowing the bending measurement to depend arbitrarily on the curvature. Moreover, in order to preserve the mechanical meaning of the boundary energy, the symmetries between the normal and geodesic curvatures of the boundary are respected, so that this energy depends on the (Frenet) curvature and torsion.

The governing equation on the interior of the surface $\Sigma$ is the Euler-Lagrange equation associated to $E[X]$, \eqref{E}, given by
\begin{equation}\label{PDE}
\Delta\dot{P}+2\dot{P}\left(2H^2-K\right)-4HP=0\,,
\end{equation}
which is a fourth-order partial differential equation in general. The other Euler-Lagrange equations associated to $E[X]$, \eqref{E}, represent three boundary conditions relating the curvatures and torsion of the boundary curve. In the general case, these boundary conditions can be understood as fixing both Dirichlet and Neumann conditions for the mean curvature in the boundary, together with an extra compatibility equation. Consequently, determining equilibria for $E[X]$, \eqref{E}, is a (strongly) overdetermined boundary value problem. 

The study of overdetermined boundary value problems began with the work of J. Serrin \cite{Serrin}, who employed the Alexandrov moving plane method to show certain symmetries of the solutions to a potential-theoretic problem.  Methods based on the strong maximum principle have also found success in this domain, but they do not apply to the case at hand since the relevant PDE is of order four. In fact, we obtain some solutions which are not axially symmetric later on in this work. Moreover, it is interesting to note that S. Eichmann proved in \cite{Eichmann} using varifold techniques that a well-posed Helfrich boundary value problem always has a (branched) Helfrich solution everywhere on the given domain, except at a finite number of branched points. This solution is continuous everywhere and smooth away from these points.

The current contribution begins in Section 2, where the Euler-Lagrange equations associated to the total energy $E[X]$, \eqref{E}, are computed. These equations characterize equilibrium configurations, and lead (through rescalings of the immersion) to a first necessary condition for the existence of equilibria.  Following this, critical domains with axial symmetry are studied in Section 3.  Note that although it is possible to obtain solutions which are not axially symmetric to this variational problem, axially symmetric surfaces are essential, appearing in physical settings as well as various conjectures regarding the minimizers of energy functionals. In particular, F. Pedit's ``Stability Conjecture" (see e.g. \cite{Pedit}) asserts that the only stable Willmore torus is the Clifford torus. We point out here that a special class of axially symmetric Willmore surfaces was studied in \cite{VM}, although it is not related to our approach. More closely related to the present problem is Conjecture 4.1 of \cite{Palmer-Pampano-2}, which states that axially symmetric critical disks attain the infimum (whenever it exists) of the Helfrich energy with elastic boundary.  In Section 3, a first integral of the Euler-Lagrange equation \eqref{PDE} on the interior of the surface is computed, and it is shown that there are three essentially different families of axially symmetric disks critical for $E[X]$, \eqref{E}, depending on the boundary circle. Some figures are also provided to illustrate the effect of the energy density $P(H)$ on the appearance of such critical surfaces.

Section 4 restricts to the special case of energy densities of the type $P(H)=\sigma(H-c_o)^p$ for $\sigma>0$, $c_o\in\mathbb{R}$ and non-negative integer $p$, i.e. p-Willmore energies with spontaneous curvature \cite{Gruber-Toda-Tran-1}.   Equilibria are studied for each possible choice of the parameter $p$, paying special attention to the ``ground state" $H\equiv H_o=c_o$ which appears when $p\geq 2$. In this case, it is shown that the boundary is a closed (and simple, if the surface is embedded) center line of a generalized Kirchhoff elastic rod.  Moreover, using an argument due to J. Nitsche \cite{Nitsche} involving the Hopf differential \cite{Hopf}, it is shown that these equilibria are planar disks when $\eta \neq 0$. Conversely, it is also shown that when $\eta=0$ the equilibrium configurations for $E[X]$, \eqref{E}, can be obtained by solving the (volume constrained) Plateau Problem, \cite{Hildebrant} (see also \cite{Heinz}). The work is concluded with the construction of some non-trivial examples of this last case.

\section{Variational Problem}

Let $\Sigma$ be a topological disk, i.e. a connected, compact surface with one boundary component, and consider an immersion of $\Sigma$ in the Euclidean 3-space $\mathbb{R}^3$,
$$X:\Sigma\rightarrow\mathbb{R}^3\,.$$
Unless otherwise stated, it is assumed that $X(\Sigma)$ is an oriented surface of class at least $\mathcal{C}^4$ with sufficiently smooth boundary, denoted $\partial\Sigma$. Moreover, $\nu$ is used to denote the unit normal vector field along $\Sigma$ as it points outside of a convex domain, and the boundary $\partial\Sigma$ is considered to be positively oriented with respect to this choice.

\begin{rem} Although $X(\Sigma)$ is only assumed to be an immersed surface, we will focus our attention on embedded surfaces since, due to the principle of non-interpenetration of matter, they are physically the most relevant ones. The differences are mentioned in the statements of the subsequent results.
\end{rem}

Note that the boundary can be represented by an arc length parameterized curve $C$. For a sufficiently smooth $C:I\rightarrow\mathbb{R}^3$, we denote the arc length parameter of $C$ by $s\in I=\left[0,L\right]$, where $L$ denotes the curve length. Using $\left(\,\right)'$ to denote the derivative with respect to the arc length, it follows that the vector field $T(s):=C'(s)$ is the unit tangent to $C$ and the (Frenet) curvature of $C$, denoted $\kappa$, is defined by $\kappa(s):=\lVert T'(s)\rVert\geq 0$.

We consider the Frenet frame $\{T,N,B\}$ along the curve $C(s)$, where $N$ and $B$ denote the unit normal and unit binormal to $C$, respectively. Since $C$ represents $\partial\Sigma$, it must be a closed curve and this frame is well defined as the curvature of $C$ can vanish only at isolated points. Moreover, since $C$ is closed, we remark that both the curvature $\kappa(s)$ and the (Frenet) torsion $\tau(s)$ represent periodic functions. The usual Frenet equations involving the curvature $\kappa$ and torsion $\tau$ of a curve $C(s)$ are then expressed as
\begin{equation}\label{feq}
\begin{pmatrix}
T\\N\\B \end{pmatrix}'=
 \begin{pmatrix}
0 &\kappa& 0\\
-\kappa&0&\tau\\
0&-\tau&0
\end{pmatrix}\begin{pmatrix}
T\\N\\B \end{pmatrix},
\end{equation}
in terms of the Frenet frame.

As suggested by S. Germain, the free energy governing the physical system associated with an elastic plate should be measured by an integral over the plate surface, \cite{Germain} with integrand consisting of a symmetric polynomial of even-degree in the principal curvatures. By a classical theorem of Newton, this is equivalent to energies of the type
$$\mathcal{E}[X]:=\int_\Sigma \mathcal{F}(H,K)\,d\Sigma\,,$$
where $\mathcal{F}$ is a polynomial in the symmetric invariants $H$ and $K$. More generally, $\mathcal{F}$ is often considered to be real analytic (see e.g. \cite{Gruber-Toda-Tran-2}) in the mean curvature $H$ and the Gaussian curvature $K$ of the immersion $X:\Sigma\rightarrow\mathbb{R}^3$.

One of the most important functionals of this type is the conformal Willmore energy (see e.g. \cite{Willmore}), which for surfaces with boundary can be expressed as (c.f. \eqref{confwillmore})
$$\mathcal{W}[X]:=\int_\Sigma\left(H^2-K\right)d\Sigma\,.$$
This case suggests that energies where the mean and Gaussian curvatures are separated and linear in the Gaussian curvature merit further investigation. There is some physical motivation for this, since the term measuring the total Gaussian curvature has been suggested as one of the major contributions to the free energy of fusion stalks in phospholipids, as well as the tendency for membranes at phase boundary to form intermediate structures (see e.g. \cite{Siegel-Kozlov}). For these reasons, this paper considers the (interior) surface energy of the immersion $X:\Sigma\rightarrow\mathbb{R}^3$ to be $\mathcal{E}[X]$ where $\mathcal{F}(H,K):=P(H)+\eta K$, for a real analytic function $P$ of one variable defined on a suitable domain and a constant $\eta\in\mathbb{R}$.

The boundary will model a (generalized) Kirchhoff elastic rod of negligible thickness, and hence it is assumed to be the corresponding center line. Recall that a generalized Kirchhoff elastic rod is a thin elastic rod with circular cross sections and uniform density, naturally straight and prismatic when unstressed, and which is held bent and twisted by external forces or momenta acting only at its ends. 

Mathematically, generalized Kirchhoff elastic rods are modeled by a curve in space, the center line, and a choice of orthonormal frame for the normal bundle to the curve, known as the material frame. Their energy is measured by a linear combination of the bending energy of the center line and a second term measuring the twisting of the rod. In the classical sense, the bending energy is the integral of the squared curvature, i.e. $\kappa^2$. In the same fashion as for the interior of the surface, it is natural to assume the bending energy is the integral of an arbitrary function depending on the curvature $\kappa$. On the other hand, twisting is accounted for using a second term, which is the integral of the squared norm of the normal derivative of any vector field in the material frame.

More precisely, denoting the material frame by $M:=\{M_1,M_2\}$ and the center line by $C$, the energy of a generalized inextensible Kirchhoff elastic rod is given by
\begin{equation}\label{Kirchhoff}
\mathcal{K}[(C,M)]:=\int_C \left(\Lambda(\kappa)+a\lVert \nabla^\perp M_1\rVert^2+b\right)ds\,,
\end{equation}
where $\Lambda$ is a real analytic function defined on an adequate domain and $a$ and $b$ are real constants. The constant $a$ represents the torsional rigidity of the rod, while $b$ can be understood as a Lagrange multiplier restricting the rod's length. Observe that classical Kirchhoff elastic rods (e.g. \cite{Langer-Singer}) usually have a constant coefficient in front of the bending energy\textemdash the so-called flexural rigidity of the rod. The present generalization assumes that this coefficient is included in the function $\Lambda$.

Adapting the computations of \cite{Langer-Singer} for the case $\Lambda(\kappa)=\kappa^2$ and varying the material frame while fixing the center line leads to the energy of the center line, which is given by
$$\mathcal{K}[C]:=\int_C \left(\Lambda(\kappa)+\varpi\tau+\beta\right)ds\,,$$
for suitable real constants $\varpi$ and $\beta$ which depend on $a$ and $b$. Note that the torsion of the choice of material frame is encoded in the energy parameter $\varpi$.

Consequently, for an immersion $X:\Sigma\rightarrow\mathbb{R}^3$, we consider the total energy
\begin{equation}\label{energy}
E[X]:=\int_\Sigma P(H)\,d\Sigma+\eta\int_\Sigma K\,d\Sigma+\oint_{\partial\Sigma}\left(\Lambda(\kappa)+\varpi\tau+\beta\right)ds\,,
\end{equation}
where $\eta$, $\varpi$ and $\beta$ are arbitrary real constants, and the functions $P$ and $\Lambda$ are real analytic functions of one variable defined in adequate domains. The energy parameters and functions are motivated by the physics. Note that in general this energy is not conformally-invariant and the regularity methods used in the literature for conformal Willmore-type energies do not apply. 

It is useful to compute the first variation of the total energy $E[X]$, \eqref{energy}. Consider arbitrary variations of the immersion $X:\Sigma\rightarrow\mathbb{R}^3$, i.e. $X+\epsilon\delta X+\mathcal{O}(\epsilon^2)$ (on the boundary $\partial\Sigma$, which is represented by a closed curve $C$, we will denote the restriction of $\delta X$ by $\delta C$). Then, for each term in the energy, we have the following variation formulas (see e.g. Appendix of \cite{Gruber} or \cite{Palmer-Pampano-2}):
\begin{itemize}
\item \textbf{The Mean Curvature Energy:}\\
Using the standard formula for the pointwise variation of the mean curvature,
$$\delta H=\frac{1}{2}\Delta\left(\nu\cdot \delta X\right)+\left(2H^2-K\right)\nu\cdot\delta X+\nabla H\cdot \delta X\,,$$
together with $\delta\left(d\Sigma\right)=\left(-2H\nu\cdot\delta X+\nabla\cdot\left[\delta X\right]^T\right)d\Sigma$, we obtain
\begin{eqnarray*}
\delta\left(\int_\Sigma P(H)\,d\Sigma\right)&=&\int_\Sigma \dot{P}(H)\,\delta H\,d\Sigma+\int_\Sigma P(H)\,\delta\left(d\Sigma\right)\\
&=&\int_\Sigma \left(\frac{1}{2}\Delta \dot{P}+\dot{P}\left[2H^2-K\right]-2H P\right)\nu\cdot\delta X\,d\Sigma\\
&+&\oint_{\partial\Sigma}\left(P n\cdot \delta C+\frac{1}{2}\dot{P}\partial_n\left[\nu\cdot\delta X\right]-\frac{1}{2}\partial_n\dot{P}\nu\cdot\delta C\right)ds\,.
\end{eqnarray*}
Here, $\dot{P}$ denotes the derivative of $P$ with respect to $H$, $n$ is the (outward pointing) conormal to the boundary defined by $n:=T\times\nu$, and $\partial_n$ represents the derivative in this direction.
\item \textbf{The Total Gaussian Curvature:}\\
Using the standard formula for the pointwise variation of the Gaussian curvature,
$$\delta K=\nabla\cdot\left(\nabla\left[\nu\cdot\delta X\right]\intprod\left[d\nu+2H\, {\rm Id}\right]\right)+2HK\nu\cdot\delta X+\nabla K\cdot \delta X\,,$$
and an argument as above involving integration by parts, it follows that
\begin{eqnarray*}
\delta\left(\int_\Sigma K\,d\Sigma\right)&=&\oint_{\partial\Sigma}\left(\nabla\left[\nu\cdot\delta X\right]\intprod\left[d\nu+2H\, {\rm Id}\right]+K\delta X\right)\cdot n\,ds\\
&=&\oint_{\partial\Sigma}\left(\left[\tau_g'\nu+Kn\right]\cdot \delta C+\kappa_n\partial_n\left[\nu\cdot \delta X\right]\right)ds\,,
\end{eqnarray*}
where $v\intprod\omega = \omega(v)$ denotes the interior product of a vector field with a differential form (i.e., contraction of the differential form with vector field) and $\kappa_n$ and $\tau_g$ are, respectively, the normal curvature and the geodesic torsion of the boundary with respect to the immersion $X$. These functions are defined by $\kappa_n:=T'\cdot\nu$ and $\tau_g:=n'\cdot \nu$, respectively.
\item \textbf{The Boundary Energy:}\\
By standard arguments (adapt, for instance, the computations of \cite{Langer-Singer}), the first variation of the boundary energy is
$$\delta\left(\oint_{\partial\Sigma}\left[\Lambda(\kappa)+\varpi\tau+\beta\right]ds\right)=\oint_{\partial\Sigma}J'\cdot \delta C\,ds\,,$$
where the vector field $J$ is defined along $\partial\Sigma$ as
\begin{equation}\label{J}
J:=\left(2\kappa\dot{\Lambda}-\Lambda-\beta\right)T+\frac{d}{ds}\left(\frac{1}{\kappa}\dot\Lambda\, T'\right)-\varpi\, T\times T'\,.
\end{equation}
Here, $\dot{\Lambda}$ represents the derivative of $\Lambda$ with respect to $\kappa$.
\end{itemize}
Combining the information above, the first variation formula for the total energy $E[X]$, \eqref{energy}, can be expressed as
\begin{eqnarray*}
\delta E[X]&=&\int_\Sigma \left(\frac{1}{2}\Delta \dot{P}+\dot{P}\left[2H^2-K\right]-2H P\right)\nu\cdot\delta X\,d\Sigma+\oint_{\partial\Sigma}\left(\frac{1}{2}\dot{P}+\eta\kappa_n\right)\partial_n\left[\nu\cdot\delta X\right]ds\\
&&+\oint_{\partial\Sigma}\left(J'+\left[\eta K+P\right]n+\left[\eta\tau_g'-\frac{1}{2}\partial_n\dot{P}\right]\nu\right)\cdot\delta C\,ds\,.
\end{eqnarray*}
As a consequence, the Euler-Lagrange equations for equilibria of the total energy $E[X]$, \eqref{energy}, are
\begin{eqnarray}
\Delta\dot{P}+2\dot{P}\left(2H^2-K\right)-4HP&=&0\,,\quad\quad\quad\text{on $\Sigma$}\,,\label{EL1}\\
\dot{P}+2\eta\kappa_n&=&0\,,\quad\quad\quad\text{on $\partial\Sigma$}\,,\label{EL2}\\
2J'\cdot\nu+2\eta\tau_g'-\partial_n\dot{P}&=&0\,,\quad\quad\quad\text{on $\partial\Sigma$}\,,\label{EL3}\\
J'\cdot n+\eta K+P&=&0\,,\quad\quad\quad\text{on $\partial\Sigma$}\,.\label{EL4}
\end{eqnarray}
The Euler-Lagrange equation \eqref{EL1} comes from considering compactly supported variations of the interior. Equations \eqref{EL2} and \eqref{EL3} can be obtained from normal variations $\delta X=\psi\nu$. Considering variations where $\psi$ vanishes on the boundary leads to \eqref{EL2}, while the vanishing of $\partial_n\psi$ on $\partial\Sigma$ yields \eqref{EL3}. Finally, the Euler-Lagrange equation \eqref{EL4} is deduced from variations that are purely tangential to the immersion.

\begin{rem}\label{expand} Using the definition of the normal curvature $\kappa_n$ and the geodesic curvature $\kappa_g$ ($\kappa_g:=T'\cdot n$) together with \eqref{J}, it is possible to expand the terms $J'\cdot \nu$ and $J'\cdot n$ in the Euler-Lagrange equations \eqref{EL3} and \eqref{EL4}, respectively, obtaining:
\begin{eqnarray*}
J'\cdot\nu&=&(G\kappa_n)''+(\kappa\dot{\Lambda}-\Lambda-\beta)\kappa_n+(G\tau_g+\varpi)(\kappa_g'-\kappa_n\tau_g)+2G'\kappa_g\tau_g+G(\kappa_g\tau_g)',\\
J'\cdot n&=&(G\kappa_g)''+(\kappa\dot{\Lambda}-\Lambda-\beta)\kappa_g-(G\tau_g+\varpi)(\kappa_n'+\kappa_g\tau_g)-2G'\kappa_n\tau_g-G(\kappa_n\tau_g)',
\end{eqnarray*}
where for simplicity we denote $G:=\dot{\Lambda}/\kappa$.
\end{rem}

A useful necessary condition for equilibria can be obtained from variations $\delta X\equiv X$ which are dilations of the reference immersion. Under a rescaling  $X\mapsto RX$ for $R>0$, the mean curvature $H$ and the (Frenet) curvature of the boundary $\kappa$ rescale like $R^{-1}$, while the line element $ds$ and the area element $d\Sigma$, rescale linearly and quadratically, respectively. Conversely, observe that the total Gaussian curvature and the total torsion of the boundary remain unchanged under such a transformation. This information leads to the following result.

\begin{prop}\label{rescaling} Let $X:\Sigma\rightarrow\mathbb{R}^3$ be a critical immersion for $E[X]$, \eqref{energy}. Then, it follows that
\begin{equation}\label{rescalingequation}
\int_\Sigma\left(2P-H\dot{P}\right)d\Sigma=\oint_{\partial\Sigma}\left(\kappa\dot{\Lambda}-\Lambda-\beta\right)ds.
\end{equation}
\end{prop}
\textit{Proof.} For a rescaled immersion $RX$, $R>0$, the energy \eqref{energy} is given by
$$E[RX]=\int_\Sigma P\left(\frac{H}{R}\right)R^2\,d\Sigma+\eta\int_\Sigma K\,d\Sigma+\oint_{\partial\Sigma}\left(\Lambda\left[\frac{\kappa}{R}\right]+\varpi\frac{\tau}{R}+\beta\right)R\,ds\,.$$
Differentiating this expression with respect to $R$ at the critical immersion $X:\Sigma\rightarrow\mathbb{R}^3$ characterized by $R=1$ then yields
$$0=\int_\Sigma \left(2P-H\dot{P}\right)d\Sigma-\oint_{\partial\Sigma}\left(\kappa\dot{\Lambda}-\Lambda-\beta\right)ds\,,$$
proving the statement. \hfill$\square$
\\

The result of Proposition \ref{rescaling} represents a nice extension of previous similar results (see e.g. \cite{Gruber-Toda-Tran-2}). Note that the present extension appears due to the consideration of boundary energies. A similar expression for the Helfrich energy, i.e. $P(H)=\sigma\left(H-c_o\right)^2$ with $\sigma>0$ and $c_o\in\mathbb{R}$, with elastic boundary can be found in \cite{Palmer-Pampano-2}.

Finally, it remains to prove the existence and uniqueness of planar critical domains under some assumptions. Recall that the contact angle $\theta\in\left[-\pi,\pi\right]$ is defined as the oriented angle between $N$, the (Frenet) normal to the boundary $C$, and $\nu$, the normal to the surface $\Sigma$. This definition allows for the expression of the normal part of the Darboux frame $\{\nu,n\}$ in terms of the normal part of the Frenet frame $\{N,B\}$. Indeed, it follows that $\nu+i n=e^{i\theta}(N+iB)$, which yields the relations
\begin{eqnarray}
\kappa_g&=&\kappa\sin\theta\,,\label{kg}\\
\kappa_n&=&\kappa\cos\theta\,,\label{kn}\\
\tau_g&=&\theta'-\tau\,.\label{taug}
\end{eqnarray}

We first consider the case where the surface is ``perfectly wetted'', i.e. $\theta\equiv\pm\pi/2$. When this happens, \eqref{kn} implies that  $\kappa_n\equiv 0$ holds along the boundary of the surface, so that the boundary is represented by an asymptotic curve. In this setting, the shape of ``perfectly wetted'' equilibria of the total energy $E[X]$, \eqref{energy}, is highly restricted.

\begin{prop}\label{Ph-} Assume that $\ddot{P}\neq 0$ on $\Sigma$ and there exists a constant mean curvature solution $H\equiv H_o$ to \eqref{EL1} given implicitly by $\dot{P}(H_o)\equiv 0$. If $X:\Sigma\rightarrow\mathbb{R}^3$ is an immersion critical for $E[X]$, \eqref{energy}, whose boundary is a circle along which $\kappa_n\equiv 0$ holds, then the surface is a planar disk. Moreover, the curvature of the boundary circle satisfies
$$\left(\kappa\dot{\Lambda}-\Lambda-\beta\right)\kappa-P(0)=0\,.$$
In particular, if $P(H)=0$ holds on $\partial\Sigma$, then $\kappa\dot{\Lambda}-\Lambda-\beta=0$ holds also.
\end{prop}
\textit{Proof.} Let $X:\Sigma\rightarrow\mathbb{R}^3$ be a critical immersion for $E[X]$, \eqref{energy}, with circular boundary and assume that $\kappa_n\equiv 0$ holds along $\partial\Sigma$. This implies from \eqref{kn} that the contact angle is $\theta\equiv\pm\pi/2$, and hence from \eqref{kg}-\eqref{taug}, $\kappa_g$ is constant along $\partial\Sigma$ and $\tau_g\equiv 0$.

With this, the Euler-Lagrange equations \eqref{EL2} and \eqref{EL3} can be simplified, obtaining
$$\dot{P}=0\,,\quad\quad\quad \ddot{P}\partial_n H=0\,,$$
along $\partial\Sigma$.

Next, since $\ddot{P}\neq 0$ and the constant mean curvature $H\equiv H_o$ given by $\dot{P}(H_o)\equiv 0$ is a solution to \eqref{EL1}, the uniqueness part of the Cauchy-Kovalevskaya Theorem implies that $H\equiv H_o$ holds on $\Sigma$, i.e. the surface has constant mean curvature. Note that this theorem assumes that the differential equation has real analytic coefficients which follows from $P$ being real analytic. 

Moreover, from Proposition 5.1 of \cite{Palmer-Pampano-2} (a constant mean curvature immersion whose boundary is a circle on which $\tau_g\equiv 0$ holds is axially symmetric) it follows that $\Sigma$ is a part of a Delaunay surface, i.e. $\Sigma$ is either a planar disk or a spherical cap. In the latter case, $\kappa_n\neq 0$ along the boundary, contradictory to assumption. Therefore, the surface is a planar disk.

Finally, since $H=K=0$ holds on $\Sigma$, equation \eqref{EL4} reduces to
$$(\kappa\dot{\Lambda}-\Lambda-\beta)\kappa-P(0)=0\,,$$
where $\kappa_g\equiv-\kappa$ (note that $\theta\equiv\pi/2$ is not possible due to the orientation of the boundary). This finishes the proof. \hfill$\square$

\begin{rem}\label{Ph} The restriction $\ddot{P}\neq 0$ is quite natural in practice. Indeed, it is satisfied by the p-Willmore energies (c.f. \cite{Gruber-Toda-Tran-1}) with spontaneous curvature given by $P(H)=\sigma(H-c_o)^p$ for any $\sigma>0$, $c_o\in\mathbb{R}$ and $p>1$. To the contrary, the cases where $\ddot{P}=0$ will be studied in detail in Sections \ref{4.1} and \ref{4.2}.
\end{rem}

\section{Axially Symmetric Critical Immersions}

Axially (or rotationally) symmetric equilibria play an essential role in the understanding of the energy $E[X]$, \eqref{energy}, as mentioned in the Introduction, so this section deals with them exclusively.

Assume that $X:\Sigma\rightarrow\mathbb{R}^3$ is an axially symmetric critical immersion. Denoting the canonical coordinates of $\mathbb{R}^3$ by $\left(x,y,z\right)$, we may assume (perhaps after a rigid motion) that the $z$-axis is the axis of rotation. By criticality, it follows that \eqref{EL1} holds on $\Sigma$. 

We denote by $C_i$, $i=1,2$ the two horizontal parallels bounding any domain $\Omega$ in $\Sigma$ and consider a variation of the immersion $X$ in the constant direction $E_3:=(0,0,1)$, i.e. $\delta X\equiv E_3$. By a computation similar to those in Section 2, the first variation formula for the mean curvature energy becomes
\begin{eqnarray*}
\delta\left(\int_\Omega P(H)\,d\Sigma\right)=\oint_{\partial\Omega}\left(P\, n_3+\frac{1}{2}\dot{P}\,\partial_n \nu_3-\frac{1}{2}\partial_n\dot{P}\,\nu_3\right)ds\,.
\end{eqnarray*}
Since $\delta X\equiv E_3$ is a constant vector field which generates a symmetry of the Lagrangian, it follows that the first variation vanishes, yielding
$$0=\oint_{\partial\Omega}\left(P\, n_3+\frac{1}{2}\dot{P}\,\partial_n\nu_3-\frac{1}{2}\partial_n\dot{P}\,\nu_3\right)ds\,.$$
For $i=1,2$, this becomes
\begin{equation}\label{fir}
r\left(P\,\partial_n z+\frac{1}{2}\dot{P}\,\partial_n \nu_3-\frac{1}{2}\partial_n\dot{P}\,\nu_3\right)=\pm a\,,
\end{equation}
for some real constant $a$. Observe that $n_3 := E_3\cdot n = \nabla z\cdot n = \partial_n z$.

Since the surface $\Sigma$ is a topological disk, we can let $C_1$ shrink to a point to conclude that $a=0$. Hence, the above expression \eqref{fir} implies the equation
\begin{equation}\label{fias}
P\,\partial_n z+\frac{1}{2}\dot{P}\,\partial_n\nu_3-\frac{1}{2}\partial_n\dot{P}\,\nu_3=0\,.
\end{equation}
Since this is done for an arbitrary subdomain $\Omega$ of $\Sigma$, equation \eqref{fias} holds everywhere on $\Sigma$ and represents a first integral of \eqref{EL1}.

\begin{rem}\label{Helfrichrem} For the case of the Helfrich energy, i.e. $P(H)=\sigma\left(H-c_o\right)^2$ with $\sigma>0$ and $c_o\in\mathbb{R}$, equation \eqref{fias} can be integrated wherever $H\neq c_o$. Up to rigid motions, this yields the equation
$$H-c_o=-\frac{\nu_3}{z}\,.$$
This expression has been obtained first in \cite{PPNew}.
\end{rem}

Note that along parallels in the surface, $\partial_n z$ and $\nu_3$ are related, respectively, to the normal curvature $\kappa_n$, and the geodesic curvature $\kappa_g$. Indeed, we have the following well known relations, which we include here for completeness.

\begin{lem} Let $X:\Sigma\rightarrow\mathbb{R}^3$ be an axially symmetric immersion parameterized, using cylindrical coordinates, by $X(\varsigma,t)=\left(r(\varsigma)\cos t, r(\varsigma)\sin t,z(\varsigma)\right)$, where $\varsigma$ is the arc length parameter of the profile curve $\gamma(\varsigma):=X(\varsigma, t_o)$. Then, the geodesic curvature $\kappa_g$ and the normal curvature $\kappa_n$ of the parallels $C(t):=X(\varsigma_o,t)$ are given, respectively, by  
\begin{eqnarray}
\kappa_g&=&\pm\frac{\nu_3}{r}\,,\label{kgas}\\
\kappa_n&=&\mp\frac{\partial_n z}{r}\,,\label{knas}
\end{eqnarray}
where $r\equiv r(\varsigma_o)$ is the radius of the parallel circle. 
\end{lem}
\textit{Proof.} The proof follows immediately from the parameterization $X(\varsigma,t)$ and the definition of the geodesic and normal curvatures, $\kappa_g:=T'\cdot n$ and $\kappa_n:=T'\cdot \nu$, respectively. \hfill$\square$
\\

Assume that $\kappa_n\equiv 0$ holds along the boundary of an axially symmetric disk critical for $E[X]$, \eqref{energy}. Provided that the conditions on the function $P$ of Proposition \ref{Ph-} are satisfied, we conclude that the surface is a planar disk and that the curvature $\kappa$ of the boundary circle satisfies
$$(\kappa\dot{\Lambda}-\Lambda-\beta)\kappa-P(0)=0\,.$$ 
In particular, if $P(H)=\sigma\left(H-c_o\right)^p$, $\sigma>0$, $c_o\in\mathbb{R}$ and $p>1$ (see Remark \ref{Ph}), it is apparent from \eqref{EL2}-\eqref{EL4} that $c_o=0$ must hold and hence $P(0)=0$. Therefore, $\kappa\dot{\Lambda}-\Lambda-\beta=0$ holds on $\partial\Sigma$.

Now, consider the case where $\kappa_g\equiv 0$ holds along the boundary of an axially symmetric disk critical for $E[X]$, \eqref{energy}, i.e. when the boundary circle is a geodesic. From \eqref{EL2}-\eqref{EL4}, it follows that the equations
\begin{eqnarray}
\dot{P}+2\eta\kappa_n&=&0\,,\label{bcgeo1}\\
2(\kappa\dot{\Lambda}-\Lambda-\beta)\kappa_n-\partial_n\dot{P}&=&0\,,\label{bcgeo2}\\
P+\eta K&=&0\,,\label{bcgeo3}
\end{eqnarray}
hold along $\partial\Sigma$. Consequently, the quantity $\kappa\dot{\Lambda}-\Lambda-\beta$ may not be zero along $\partial\Sigma$.

In the rest of the section, we are going to obtain examples of axially symmetric disks critical for the total energy $E[X]$, \eqref{energy}. For that, we will solve equation \eqref{fias} and use the boundary conditions described above.

Let $\gamma(\varsigma)=\left(r(\varsigma),z(\varsigma)\right)$ be the arc length parameterized profile curve of an axially symmetric disk satisfying \eqref{fias}. Since $\gamma(\varsigma)$ is parameterized by arc length, there is a function $\varphi(\varsigma)$ such that
\begin{eqnarray}
r'(\varsigma)&=&\cos\varphi(\varsigma)\,,\label{dif1}\\
z'(\varsigma)&=&\sin\varphi(\varsigma)\,.\label{dif2}
\end{eqnarray}
The function $\varphi(\varsigma)$ represents the angle between the positive part of the $r$-axis and the tangent vector to $\gamma(\varsigma)$. With this notation, equation \eqref{fias} can be written as
\begin{equation}\label{dif3}
P(H)\sin\varphi(\varsigma)\mp\frac{1}{2}\dot{P}(H)\varphi'(\varsigma)\sin\varphi(\varsigma)\mp\frac{1}{2}\left(\dot{P}(H)\right)_\varsigma\cos\varphi(\varsigma)=0\,,
\end{equation}
where $\left(\,\right)'\equiv\left(\,\right)_\varsigma$ represents the derivative with respect to $\varsigma$ and
\begin{equation}\label{H}
H\equiv H(\varsigma)=\pm\frac{\sin\varphi(\varsigma)}{2r(\varsigma)}+\frac{\varphi'(\varsigma)}{2}\,.
\end{equation}
Observe that, if $\dot{P}'(H)=\ddot{P}H'=0$ holds on $\Sigma$, then either $H\equiv H_o$ is constant or $\ddot{P}(H)\equiv 0$. In the former, the surface is a part of a Delaunay surface, i.e. either a planar disk or a spherical cap; while, in the latter, the surface energy is either the surface area (see Section \ref{4.1}) or the $1$-Willmore energy with spontaneous curvature (see Section \ref{4.2}).

Since these two cases will be analyzed in detail later, now we will focus on the case $\ddot{P}H'\neq 0$. In this case, the differential equation \eqref{dif3} is a second order differential equation, and the profile curve $\gamma(\varsigma)$ satisfies the system of differential equations \eqref{dif1}-\eqref{dif3}.

We now impose the initial conditions corresponding to $\varsigma=0$. Since the surface is a regular topological disk, the profile curve $\gamma(\varsigma)$ must cut the $z$-axis perpendicularly, so that $r(0)=0$ and $\varphi(0)=0$. Moreover, by the invariance of these surfaces under vertical translations, we may fix the initial height to be $z(0)=1$, while keeping $\varphi'(0)$ as a parameter.

For the boundary conditions ($\varsigma=L$ where $L$ is the length of $\gamma$), we will assume that $P+\eta K=0$ holds along the boundary. Then, two cases must be distinguished depending on whether the boundary circle is a geodesic or not. In the former, i.e. when the boundary is a geodesic circle, then \eqref{bcgeo1}-\eqref{bcgeo3} must hold. In the latter, $\kappa_g\neq 0$ and it may also be assumed that $\kappa_n\neq 0$, since this case is mainly covered in Proposition \ref{Ph-}. Therefore, $\kappa\dot{\Lambda}-\Lambda-\beta=0$; $\dot{P}+2\eta\kappa_n=0$; $P+\eta K=0$ and $\partial_n\dot{P}\equiv\dot{P}_\varsigma=0$ must hold along $\partial\Sigma$, which can be simplified to the following three equations:
\begin{eqnarray}
\kappa\dot{\Lambda}-\Lambda-\beta&=&0\,,\label{bcnogeo1}\\
\dot{P}+2\eta\kappa_n&=&0\,,\label{bcnogeo2}\\
2P-(2H-\kappa_n)\dot{P}&=&0\,.\label{bcnogeo3}
\end{eqnarray}

\begin{rem} If the boundary energy is measured by the total curvature, i.e. $\Lambda(\kappa)=\alpha\kappa-\beta$ for some real constant $\alpha$, then the equation $\kappa\dot{\Lambda}-\Lambda-\beta=0$ is an identity. Recall that the total curvature is constant on a regular homotopy class of planar curves. For the rest of the cases, $\kappa\dot{\Lambda}-\Lambda-\beta=0$ represents an algebraic equation which determines the radius of the boundary circle.
\end{rem}

\begin{figure}[H]
\centering
\begin{subfigure}[b]{0.3\linewidth}
\includegraphics[width=\linewidth]{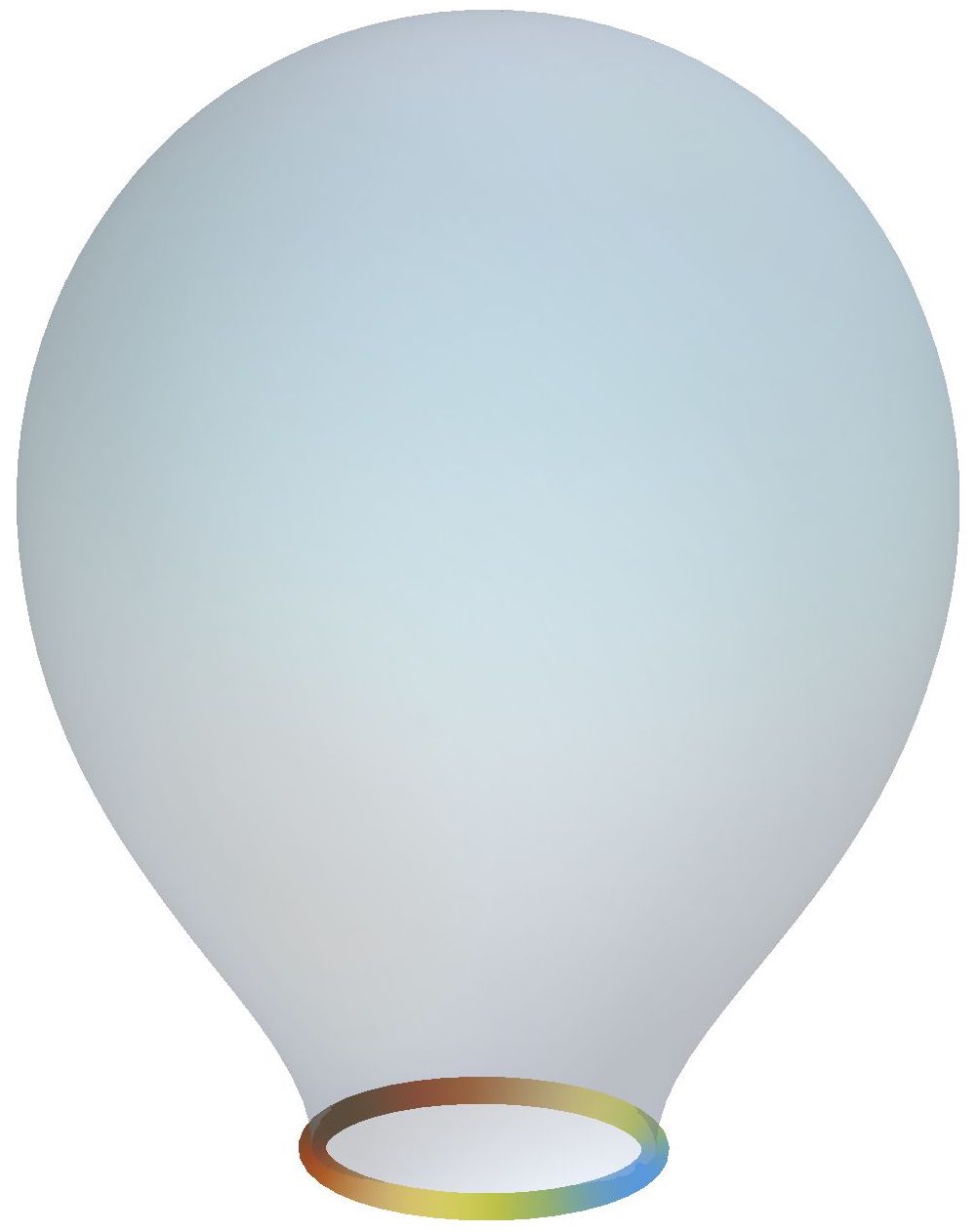}
\caption{$P(H)=e^{H^2}$}
\end{subfigure}
\quad
\begin{subfigure}[b]{0.3\linewidth}
\includegraphics[width=\linewidth]{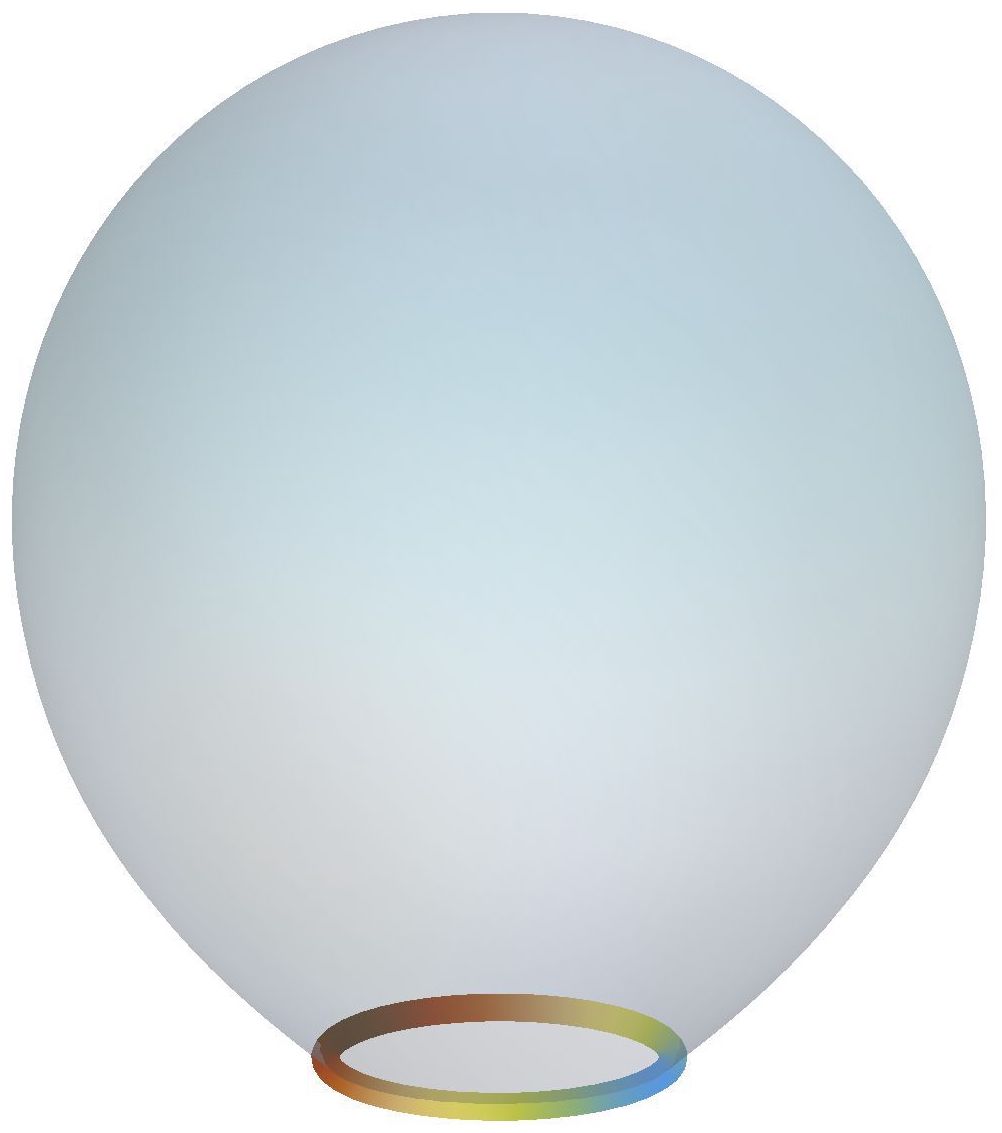}
\caption{$P(H)=(H+2)^2$}
\end{subfigure}
\quad
\begin{subfigure}[b]{0.3\linewidth}
\includegraphics[width=\linewidth]{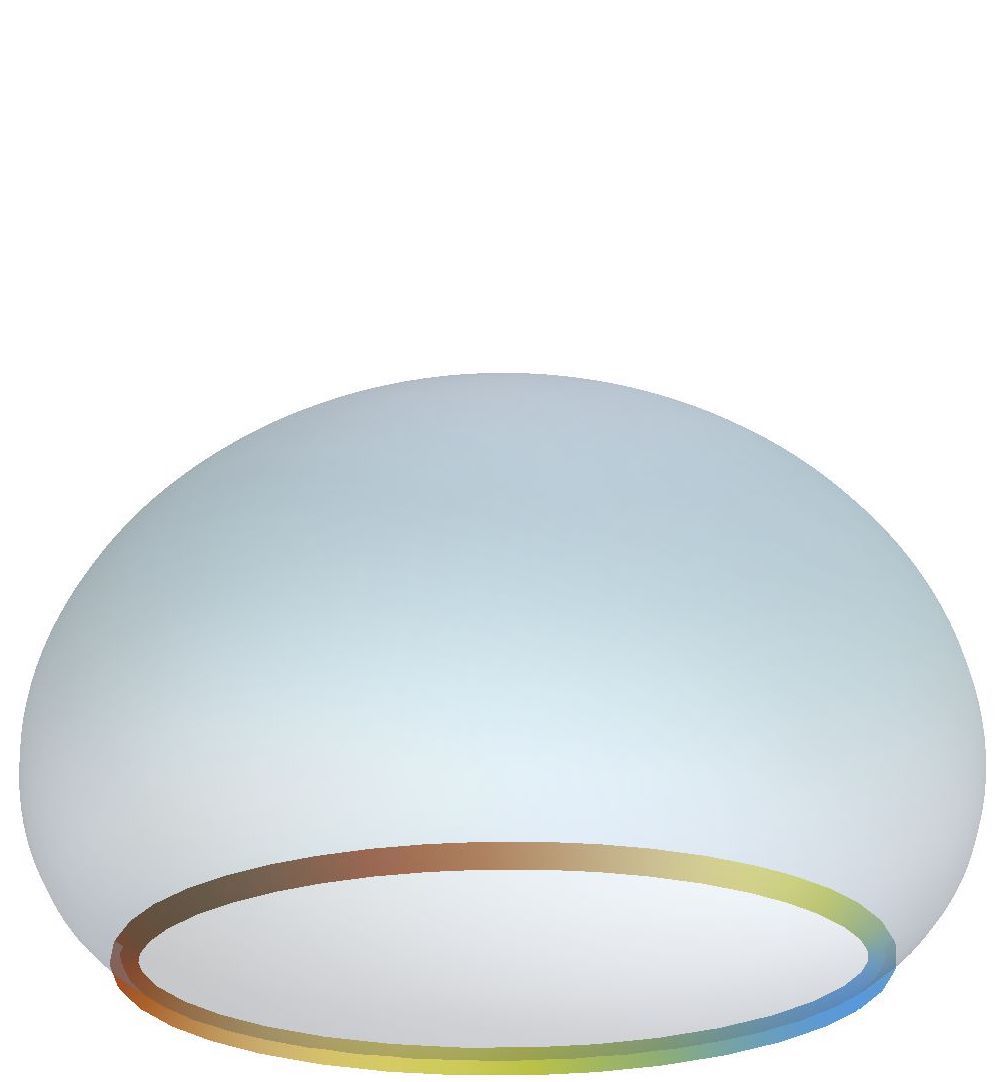}
\caption{$P(H)=\log H^2$}
\end{subfigure}
\caption{Axially symmetric critical disks for the energy $E[X]$, \eqref{energy}, for different surface energies and suitable $\Lambda(\kappa)$, $\eta$ and $\beta$. In case (a), the boundary is a geodesic circle; while for (b) and (c) the boundary is a non-geodesic suitable parallel. Case (b) was first obtained in \cite{PPNew}.}
\label{bal}
\end{figure}

In Figure \ref{bal}, we show three critical embedded domains for the energy $E[X]$, \eqref{energy}, where the surface energy is determined by different real analytic functions $P(H)$, defined on suitable domains. These critical disks are obtained by solving the system of differential equations \eqref{dif1}-\eqref{dif3} with the initial conditions described above. In case (a), the boundary circle is a geodesic parallel and along it \eqref{bcgeo1}-\eqref{bcgeo3} are satisfied for suitable energy parameters $\eta$ and $\beta$. On the other hand, for cases (b) and (c), the boundary circles are not geodesic and along them \eqref{bcnogeo1}-\eqref{bcnogeo3} are satisfied for suitable energy parameters $\eta$ and $\beta$. We point out here that for the Helfrich energy, i.e. $P(H)=\sigma (H-c_o)^2$ (see Figure \ref{bal}, (b)), more axially symmetric critical disks, obtained using a similar approach and the integral equation of Remark \ref{Helfrichrem}, can be found in \cite{PPNew}.

\begin{rem} There exists another case when $P+\eta K\neq 0$ along the boundary.
\end{rem}

\section{Generalized Willmore Surface Energies}

As mentioned in Section 2, S. Germain proposed to consider surface energies where the integrand is a symmetric (and even degree) polynomial in the principal curvatures \cite{Germain}. Interesting examples of these integrands lead to generalized Willmore energies, also called $p$-Willmore energies with spontaneous curvature (c.f. \cite{Gruber-Toda-Tran-1}). 

This section will specialize the previous findings to these particular cases and study different properties related with the ``ground state'' $H\equiv H_o=c_o$.

Consider the family of functions $P(H)=\sigma\left(H-c_o\right)^p$ for $\sigma>0$, $c_o\in\mathbb{R}$ and a non-negative integer $p$. The parameter $c_o$ represents the spontaneous curvature and, consequently, the (interior) surface energy corresponding to the integral of $P(H)$ is a p-Willmore type energy with spontaneous curvature. This yields the following particular case of the energy $E[X]$, \eqref{energy},
\begin{equation}\label{F}
F[X]:=\sigma\int_\Sigma\left(H-c_o\right)^p\,d\Sigma+\eta\int_\Sigma K\,d\Sigma+\oint_{\partial\Sigma}\left(\Lambda(\kappa)+\varpi\tau+\beta\right)ds\,.
\end{equation}
In this context, the parameter $\sigma>0$ is called the bending rigidity of the surface.

\begin{rem} We point out here a first difference between the cases $p$ even and $p$ odd. While for $p$ even $P(H)=\sigma\left(H-c_o\right)^p\geq 0$ holds, for $p$ odd the integral of this quantity may be negative (c.f. \cite{Gruber}), resulting, for suitable energy parameters $\eta$, $\varpi$, $\beta$ and the function $\Lambda(\kappa)$, in the undesired scenario of the energy $F[X]$, \eqref{F}, being negative. However, for our purposes of analyzing equilibria we will not worry about it.
\end{rem}

The Euler-Lagrange equations associated to the energy $F[X]$, \eqref{F}, can be obtained by specializing equations \eqref{EL1}-\eqref{EL4}. Indeed, the respective equations are given by
\begin{eqnarray}
\Delta\left(p\left[H-c_o\right]^{p-1}\right)+2\left(H-c_o\right)^{p-1}\left(2[p-1]H^2-pK+2c_oH\right)&=&0\,,\label{ELF1}\\
\sigma p\left(H-c_o\right)^{p-1}+2\eta\kappa_n&=&0\,,\label{ELF2}\\
2J'\cdot \nu-\sigma p(p-1)\left(H-c_o\right)^{p-2}\partial_nH+2\eta\tau_g'&=&0\,,\label{ELF3}\\
J'\cdot n+\sigma\left(H-c_o\right)^p+\eta K&=&0\,,\label{ELF4}
\end{eqnarray}
where the vector field $J$ is given in \eqref{J} (see also Remark \ref{expand} for the expanded quantities $J'\cdot n$ and $J'\cdot \nu$).

First considered will be planar critical domains. In these cases $H=K=0$ holds on the surface and $\kappa_n=\tau_g=0$ holds along the boundary. Therefore, the Euler-Lagrange equations \eqref{ELF1}-\eqref{ELF4} simplify to
$$p\left(-c_o\right)^{p-1}=0\,,\quad\quad\quad J'\cdot \nu=0\,,\quad\quad\quad J'\cdot n+\sigma\left(-c_o\right)^p=0\,.$$
It is then clear that if $p\neq 0$ (the case $p=0$ will be studied separately, see Section 4.1), $p>1$ and $H=c_o=0$ holds. That is, there is no spontaneous curvature. Moreover, it follows from Remark \ref{expand}, that
\begin{eqnarray}
\varpi\kappa'&=&0\,,\label{one}\\
\dot{\Lambda}''+\kappa^2\dot{\Lambda}-\kappa\left(\Lambda+\beta\right)&=&0\,.\label{two}
\end{eqnarray}
From equation \eqref{one}, $\kappa$ is constant if $\varpi\neq 0$, and hence the boundary is a circle. On the other hand, equation \eqref{two} is the Euler-Lagrange equation associated to the curvature energy
$$\mathbf{\Theta}[C]:=\int_C\left(\Lambda(\kappa)+\beta\right)ds\,,$$
acting on the space of planar curves. The following result summarizes this information.

\begin{prop}\label{planar} A planar disk is critical for $F[X]$, \eqref{F}, with $p\neq 0$, if and only if $p>1$, $c_o=0$ and the boundary is a closed curve contained in a plane critical for
\begin{equation}\label{new}
\mathbf{\Theta}[C]:=\int_C\left(\Lambda(\kappa)+\beta\right)ds\,,
\end{equation}
satisfying $\varpi\kappa'=0$. Moreover, if the surface is embedded, then the boundary is simple.
\end{prop}

If $\varpi\neq 0$ then, as mentioned above, the boundary is a circle whose curvature $\kappa$ satisfies the equation
$$\kappa\dot{\Lambda}-\Lambda-\beta=0\,.$$
On the contrary, if $\varpi=0$, then (apart from the previous circles) any planar closed curve critical for $\mathbf{\Theta}[C]$, \eqref{new}, may be the boundary of a planar critical domain. However, existence of planar closed critical curves for $\mathbf{\Theta}[C]$, \eqref{new}, is highly restricted. First of all, the curvature of the critical curve, i.e. a solution of \eqref{two}, must be periodic. Moreover, the following closure condition must also be satisfied:
$$\int_0^\rho\left(\kappa\dot{\Lambda}-\Lambda-\beta\right)ds=0\,,$$
where $\rho$ is the period of the curvature. This condition may be obtained solely though analyzing critical curves for $\mathbf{\Theta}[C]$, \eqref{new}. However, it is proven here using a different approach.

\begin{prop}\label{intcon} Let $X:\Sigma\rightarrow\mathbb{R}^3$ be a planar disk critical for $F[X]$, \eqref{F}, with $p\neq 0$. Then, 
$$\int_0^\rho \left(\kappa\dot{\Lambda}-\Lambda-\beta\right)ds=0\,,$$
where $\rho$ is the period of the curvature of the boundary curve. (If the boundary is a circle, $\rho$ can be assumed to be any real number.)
\end{prop}
\textit{Proof.} From Proposition \ref{planar}, for planar critical disks, $H=c_o=0$ holds on $\Sigma$. Therefore, applying the rescaling condition obtained in Proposition \ref{rescaling}, we have
$$\oint_{\partial\Sigma} \left(\kappa\dot{\Lambda}-\Lambda-\beta\right)ds=\int_\Sigma \left(2P-H\dot{P}\right)d\Sigma=\sigma\int_\Sigma\left(H-c_o\right)^{p-1}\left(2\left[H-c_o\right]-pH\right)d\Sigma=0\,.$$
Finally, notice that since the boundary curve is closed, its curvature $\kappa$ must be periodic and so is $\kappa\dot{\Lambda}-\Lambda-\beta$. Moreover, these functions share the same period $\rho$. Thus, the integral along the entire boundary $\partial\Sigma$ is a natural multiple $n\rho$ of the integral over one period, where $n\in\mathbb{N}$ is the smallest positive integer needed to close the curve. \hfill$\square$
\\

Moreover, it is also interesting to consider consequences of the same rescaling condition when the boundary energy is arranged to vanish. Recall that at a critical immersion of $\Sigma$, Proposition \ref{rescaling} guarantees the flux formula
\begin{equation}\label{flux}
0=\int_\Sigma \left(H- c_o\right)^{p-1}\left([2-p]H - 2c_o\right)d\Sigma\,,
\end{equation}
provided the boundary integral is arranged to vanish. An elementary analysis of this condition yields the following rigidity consequences.

\begin{prop} Suppose that $X:\Sigma\rightarrow\mathbb{R}^3$ is a critical immersion of $\Sigma$ for the energy $F[X]$, \eqref{F}, and 
$$\oint_{\partial\Sigma} \left(\kappa\dot{\Lambda}-\Lambda-\beta\right)ds = 0\,.$$
Then, the following hold:
\begin{enumerate}[(i)]
\item If $p=1$, 
$$\mathcal{W}^1[X]:=\int_\Sigma \left(H-c_o\right)d\Sigma=c_o\,\mathcal{A}[X]\,,$$ 
where $\mathcal{A}[X]$ denotes the surface area. 
\item If $p=2$, then either $c_o= 0$ or $\mathcal{W}^1[X] = 0$.
\item If $p>2$ is even and $\lvert H\rvert\geq \lvert c_o\rvert$ on $\Sigma$, then $\Sigma$ has constant mean curvature.
\item If $p>2$ is odd and $\lvert H\rvert \leq 2\lvert c_o\rvert/(p-2)$ on $\Sigma$, then $\Sigma$ has constant mean curvature. 
\end{enumerate}
\end{prop}
\textit{Proof.} Suppose $p=1$.  Then, the hypotheses along with Proposition \ref{rescaling} imply the equality
$$0 = \int_\Sigma \left(H - 2c_o\right)d\Sigma\,,$$
which immediately establishes the first consequence. Moreover, when $p=2$, the relevant flux formula, \eqref{flux}, becomes 
$$0 = -2c_o\int_\Sigma \left(H - c_o\right)d\Sigma = -2c_o\,\mathcal{W}^1[X]\,,$$
so either $c_o = 0$, in which case the surface energy is the (conformally invariant, and thus, rescaling invariant) Willmore energy, or $\mathcal{W}^1[X] = 0$. Note that the latter conclusion implies that the total mean curvature is a multiple of the surface area, the multiple being the spontaneous curvature.

For the remainder, suppose that $p>2$. We first prove the result for $c_o\geq 0$. Since $p>2$, it follows that whenever $H\geq c_o$ we have the inequalities
$$(2-p)H-2c_o\leq (2-p)c_o-2c_o= -pc_o \leq 0\,.$$
In particular, from this and the flux formula \eqref{flux}, when $p>2$,
$$0 = \int_\Sigma \left(H-c_o\right)^{p-1}\left([2-p]H-2c_o\right)d\Sigma\leq 0\,,$$
so that $H\equiv H_o=2c_o/(2-p)$ holds on $\Sigma$, since, due to our choice of orientation, the constant mean curvature of a surface is non positive. Similarly, whenever $H \leq -c_o$, there is the inequality
$$(2-p)H-2c_o\geq (p-4)c_o\,.$$
Therefore, from \eqref{flux}, when $p>2$ is even, it follows that
$$0 = \int_\Sigma \left(H-c_o\right)^{p-1}\left([2-p]H - 2c_o\right)d\Sigma \geq 0\,,$$
so that $H\equiv H_o= 2c_o / (2-p)$ on $\Sigma$. 

Now, for $c_o<0$, a similar argument shows that when $p> 2$ is even, $H\equiv H_o=c_o$ holds on $\Sigma$. This establishes the third statement.

Finally, assume $c_o\geq 0$. If $p>2$ is odd and $H\geq 2c_o/(2-p)$, then 
$$0 = \int_\Sigma \left(H- c_o\right)^{p-1}\left([2-p]H - 2c_o\right)d\Sigma  \leq 0\,,$$
so $H \equiv H_o=2c_o/(2-p)\leq 0$ holds and $\Sigma$ has constant mean curvature. Similarly, if $c_o<0$, we conclude that $H\equiv H_o=c_o< 0$ on $\Sigma$. This establishes the final statement. \hfill$\square$
\\

We consider in what follows different interesting cases:

\subsection{Area Functional ($p=0$)}\label{4.1}

In this section we will fix $p=0$ so that the energy $F[X]$, \eqref{F}, reduces to an extended Kirchhoff-Plateau energy with elastic modulus (see \cite{Gruber-Pampano-Toda} and references therein). In particular, the interior surface energy is the classical surface area. In this setting, the Euler-Lagrange equations \eqref{ELF1}-\eqref{ELF4} read
\begin{eqnarray}
H&\equiv&0\,,\quad\quad\quad\text{on $\Sigma$}\,,\label{EL01}\\
\eta\kappa_n&=&0\,,\quad\quad\quad\text{on $\partial\Sigma$}\,,\label{EL02}\\
J'\cdot \nu+\eta\tau_g'&=&0\,,\quad\quad\quad\text{on $\partial\Sigma$}\,,\label{EL03}\\
J'\cdot n-\eta\tau_g^2+\sigma&=&0\,,\quad\quad\quad\text{on $\partial\Sigma$}\,,\label{EL04}
\end{eqnarray}
where $J$ is defined in \eqref{J}, and the last equation used the definition of the Gaussian curvature $K$ along $\partial\Sigma$, i.e. $K=\kappa_n\left(2H-\kappa_n\right)-\tau_g^2$, together with \eqref{EL01} and \eqref{EL02}.

Observe that the Euler-Lagrange equation \eqref{EL01} implies that equilibria must be minimal, which can be understood as a variational characterization of the surface alone. On the other hand, the boundary conditions \eqref{EL03} and \eqref{EL04} can be rewritten in vector form as
\begin{equation}\label{halfint}
\left(J+\eta\tau_g\nu\right)'+\sigma n=0\,.
\end{equation}
After integrating along the boundary $C$, it follows that
$$\oint_C n\, ds=0\,,$$
must hold.

\begin{rem} This has an interesting connection to conjugate minimal surfaces. Indeed, let $Y$ be the conjugate minimal immersion to X (which is globally defined since $\Sigma$ is a topological disk). The Cauchy-Riemann equations imply that $Y'=n$, and hence \eqref{halfint} can be integrated once to obtain the third-order conservation law:
\begin{equation}\label{int}
J+\eta\tau_g\nu+\sigma Y\equiv A\,,
\end{equation}
for a constant vector $A\in\mathbb{R}^3$ (c.f. \cite{Palmer-Pampano-1}).
\end{rem}

Moreover, when $\eta\neq 0$, the Euler-Lagrange equation \eqref{EL02} determines how the minimal surface $\Sigma$ meets its boundary $\partial\Sigma$, which gives a restriction on the contact angle $\theta$ from before. Indeed, it is a consequence of $\kappa_n=\kappa\cos\theta$ and \eqref{EL02} that $\theta\equiv \pm\pi/2$, i.e. the surface is ``perfectly wetted". In fact, these conditions imply that the surface is planar.

\begin{thm}\label{t} Let $\eta\neq 0$ and $X:\Sigma\rightarrow\mathbb{R}^3$ be an immersion of a disk type surface critical for $F[X]$, \eqref{F}, with $p=0$. Then, the surface is a planar domain bounded by a (planar) closed curve whose curvature $\kappa$ satisfies $\varpi\kappa'=0$ and
\begin{equation}\label{areacons}
\dot{\Lambda}''+\kappa^2\dot{\Lambda}-\kappa\left(\Lambda+\beta\right)\pm\sigma=0\,.
\end{equation} 
Moreover, if the surface is embedded, then the boundary is simple.
\end{thm}
\textit{Proof.} Let $\Sigma$ be a topological disk and consider the immersion $X:\Sigma\rightarrow\mathbb{R}^3$  critical for $F[X]$, \eqref{F}, with $p=0$. By definition, this means the Euler-Lagrange equations \eqref{EL01}-\eqref{EL04} hold for $X$. From \eqref{EL01}, the surface is minimal. Moreover, since $\eta\neq 0$ by hypothesis, it follows from \eqref{EL02} that $\kappa_n\equiv 0$ on $\partial\Sigma$.

Next, a modification of an argument due to Nitsche (see e.g. \cite{Nitsche}) based on the holomorphicity of the Hopf differential (\cite{Hopf}) can be applied to conclude that a constant mean curvature (in particular, minimal) surface of disk type bounded by an asymptotic curve is necessarily a planar domain. For details, see Theorem 4.1 of \cite{Palmer-Pampano-2}, or also Theorem 4.1 of \cite{Gruber-Pampano-Toda}. From this, it follows that $K=H=0$ hold on $\Sigma$ and $\kappa_n=\tau_g=0$ hold on $\partial\Sigma$.

Finally, using this information in the boundary conditions \eqref{EL03} and \eqref{EL04} leads to
\begin{eqnarray*}
\varpi\kappa'&=&0\,,\\
\dot{\Lambda}''+\kappa^2\dot{\Lambda}-\kappa\left(\Lambda+\beta\right)\pm\sigma&=&0\,,
\end{eqnarray*}
respectively, proving the statement. \hfill$\square$
\\

Due to this result, if $\varpi\neq 0$ then $\kappa\equiv\kappa_o$ is constant and the boundary curve is a circle. Moreover, equation \eqref{areacons} is satisfied for the boundary of a critical disk independently of the value  $\varpi$. Closed planar curves whose curvature satisfy \eqref{areacons} can be characterized as critical curves for the energy
\begin{equation}\label{energycurves1}
\mathbf{\Theta}[C]:=\int_C\left(\Lambda(\kappa)+\beta\right)ds\,,
\end{equation}
under variations preserving the enclosed algebraic area. This restriction is included via a Lagrange multiplier depending on $\sigma>0$. Consequently, boundary curves in this setting represent area-constrained planar critical curves for $\mathbf{\Theta}[C]$, \eqref{energycurves1}. For a survey of area-constrained planar elastic curves (i.e. $\Lambda(\kappa)=\kappa^2$), see \cite{35} and references therein.

\subsection{Total Mean Curvature ($p=1$)}\label{4.2}

We will focus now on the case $p=1$ for the energy $F[X]$, \eqref{F}. The surface energy in this case becomes a multiple of the $1$-Willmore energy with spontaneous curvature,
$$\mathcal{W}^1[X] :=\int_\Sigma \left(H - c_o\right)d\Sigma\,.$$
This energy can be understood as the total mean curvature with area constraint. Of course, this reduces to the usual total mean curvature functional when $c_o=0$.

From equations \eqref{ELF1}-\eqref{ELF4} we have that the Euler-Lagrange equations characterizing equilibria for the case $p=1$ in $F[X]$, \eqref{F}, are given by
\begin{eqnarray}
K-2c_o H&\equiv&0\,,\quad\quad\quad\text{on $\Sigma$}\,,\label{EL11}\\
2\eta\kappa_n+\sigma&=&0\,,\quad\quad\quad\text{on $\partial\Sigma$}\,,\label{EL12}\\
J'\cdot \nu+\eta\tau_g'&=&0\,,\quad\quad\quad\text{on $\partial\Sigma$}\,,\label{EL13}\\
J'\cdot n+\eta K+\sigma\left(H-c_o\right)&=&0\,,\quad\quad\quad\text{on $\partial\Sigma$}\,,\label{EL14}
\end{eqnarray}
where $J$ is given in \eqref{J}. Notice that equation \eqref{EL11} implies that our surface is a linear Weingarten surface (c.f. \cite{Weingarten}). Since $\sigma>0$, it follows from equation \eqref{EL12} that $\eta\neq 0$ must hold in order equilibria to exist, and hence $\kappa_n$ is constant along the boundary. This suggests a search for axially symmetric critical immersions. A first result in this direction is the following.

\begin{prop} Let $X:\Sigma\rightarrow\mathbb{R}^3$ be an axially symmetric critical immersion for $F[X]$, \eqref{F}, with $p=1$. Then, $\sigma=-4\eta c_o>0$ and $\kappa_n=2c_o$ along $\partial\Sigma$. Moreover, the curvature $\kappa$ of the boundary curve satisfies $\kappa\dot{\Lambda}-\Lambda-\beta=0$.
\end{prop}
\textit{Proof.} For axially symmetric immersions, recall that $\kappa_g$ and $\kappa_n$ are constant along the boundary $\partial\Sigma$ while $\tau_g\equiv 0$ holds there. Moreover, the Euler-Lagrange equations \eqref{EL11}-\eqref{EL14} must hold since the surface is critical.

Equation \eqref{EL12} directly implies that $\eta\neq 0$ and $\kappa_n=-\sigma/(2\eta)\neq 0$ along $\partial\Sigma$. At the same time, equations \eqref{EL13} and \eqref{EL14} reduce, respectively, to
\begin{eqnarray*}
\kappa\dot{\Lambda}-\Lambda-\beta&=&0\,,\\
\eta K+\sigma\left(H-c_o\right)&=&0\,.
\end{eqnarray*}
Finally, combining \eqref{EL11} and $\kappa_n=-\sigma/(2\eta)$ with the definition of $K$ along $\partial\Sigma$, it follows from the second equation above that $\sigma=-4\eta c_o$ holds. \hfill$\square$
\\

In order to find specific examples, we need to search for axially symmetric surfaces satisfying $K\equiv 2c_o H$ on $\Sigma$ and such that $\kappa_n=2c_o\neq 0$ along $\partial\Sigma$. Equivalently, this means the boundary must satisfy $H=2c_o\neq 0$. 

Trivial examples of axially symmetric surfaces satisfying $K\equiv 2c_oH$ are planar disks and spherical caps. In the former, $H\equiv 0$ holds on $\Sigma$, so it is impossible to find a boundary circle along which $H=2c_o\neq 0$ holds. However, in a sphere of curvature $K\equiv 4 c_o^2$ any spherical cap satisfies above conditions. Therefore, for suitable energy parameters and radius of the boundary circle, these spherical caps are critical for $F[X]$, \eqref{F}, with $p=1$ (c.f. Proposition \ref{caps}).

We now study other possible axially symmetric solutions of above conditions. The profile curves of axially symmetric surfaces satisfying $K\equiv 2c_o H$ were characterized in \cite{Pampano} as planar critical curves for 
\begin{equation}\label{curvpr}
\mathbf{\Theta}[\gamma]:=\int_\gamma \sqrt{\varepsilon\left(\chi^2+2c_o\chi\right)}\,d\varsigma\,,
\end{equation}
where $\chi$ denotes the curvature of the profile curve $\gamma$, $\varepsilon$ represents the sign of $\chi^2+2c_o\chi$, and $\varsigma$ is its arc length parameter. The first integral of the Euler-Lagrange equation associated to $\mathbf{\Theta}[\gamma]$, \eqref{curvpr}, can be expressed as (for details, see \cite{Pampano})
\begin{equation}\label{fichi}
\left(\chi'\right)^2=\frac{\left(\chi^2+2c_o\chi\right)^2}{c_o^4}\left(\left[\varepsilon d-c_o^2\right]\chi^2+2\varepsilon c_o d\chi\right),
\end{equation}
for a positive real constant $d$.

Moreover, by classical arguments on the theory of curvature energy functionals, we have that these surfaces are binormal evolution surfaces (\cite{Pampano}) and, hence, they can be parameterized as
$$X(\varsigma,t)=\frac{\varepsilon}{\sqrt{d}}\left(\frac{\chi+c_o}{\sqrt{\varepsilon\left(\chi^2+2c_o\chi\right)}}\cos t,\frac{\chi+c_o}{\sqrt{\varepsilon\left(\chi^2+2c_o\chi\right)}}\sin t, -c_o\int \frac{\chi}{\sqrt{\varepsilon\left(\chi^2+2c_o\chi\right)}}\,d\varsigma\right).$$
In this parameterization, the normal curvature of parallels, \eqref{knas}, is expressed as
$$\kappa_n=-c_o\frac{\chi}{\chi+c_o}\,,$$
and hence $\kappa_n=2c_o$ if and only if $\chi(\varsigma)=-2c_o/3$. Since we are looking for topological disks, our profile curve must cut the axis of rotation, i.e. from above parameterization $\chi(\varsigma)=-c_o$. We point out here that a change of orientation on our profile curve carries a change of sign in the energy parameter $c_o$, so we may assume that $c_o>0$ holds. With this assumption, the curvature of our profile curve varies in the domain $\chi\in\left[-c_o,-2c_o/3\right]$. Moreover, for these values of $\chi$, $\chi^2+2c_o\chi$ is negative and so $\varepsilon=-1$. Note that, as a consequence of \eqref{fichi}, $d\geq c_o^2$ must hold.

\begin{figure}[H]
\makebox[\textwidth][c]{
\begin{subfigure}[b]{0.35\linewidth}
\includegraphics[width=\linewidth]{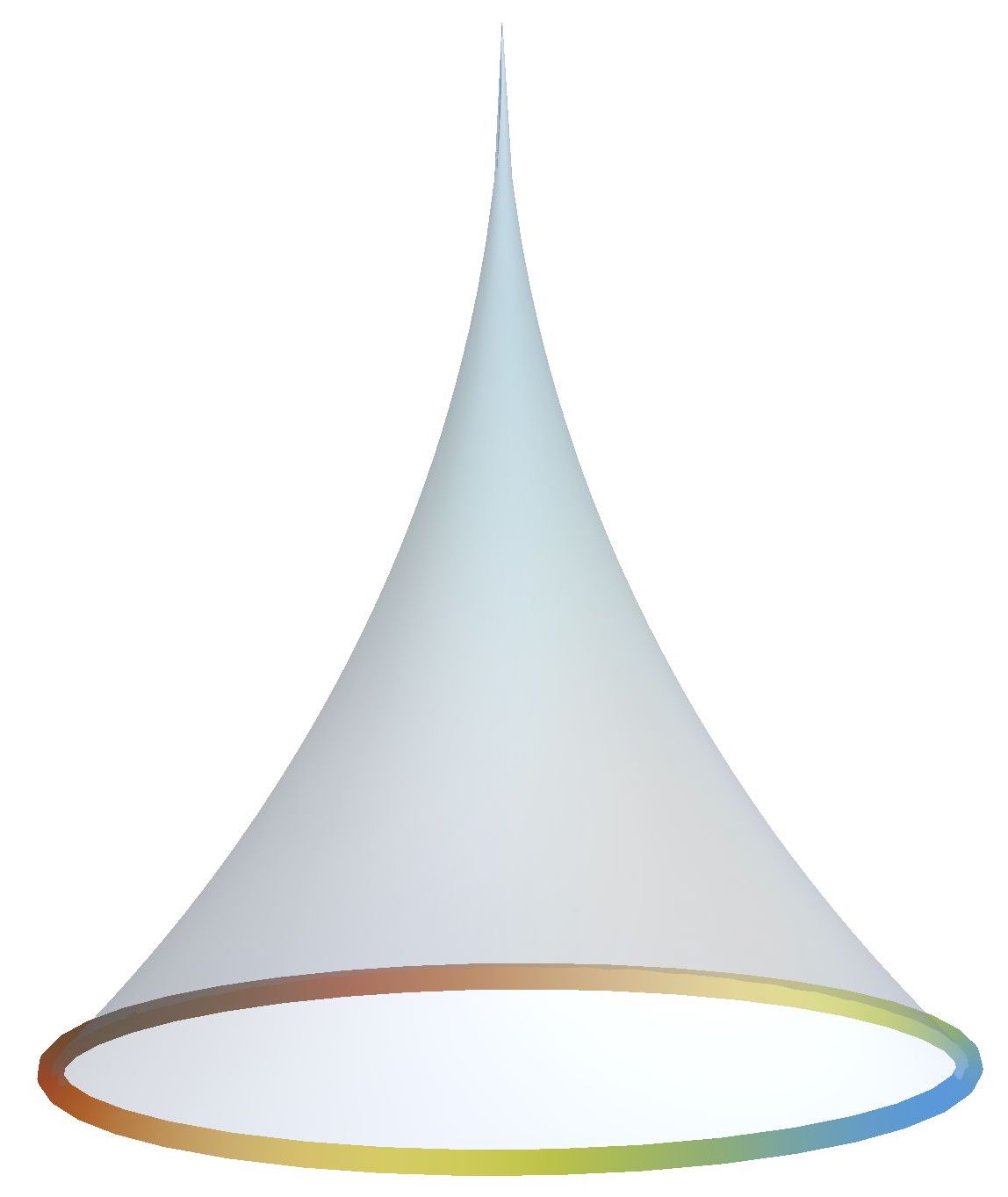}
\end{subfigure}
\begin{subfigure}[b]{0.35\linewidth}
\includegraphics[width=\linewidth]{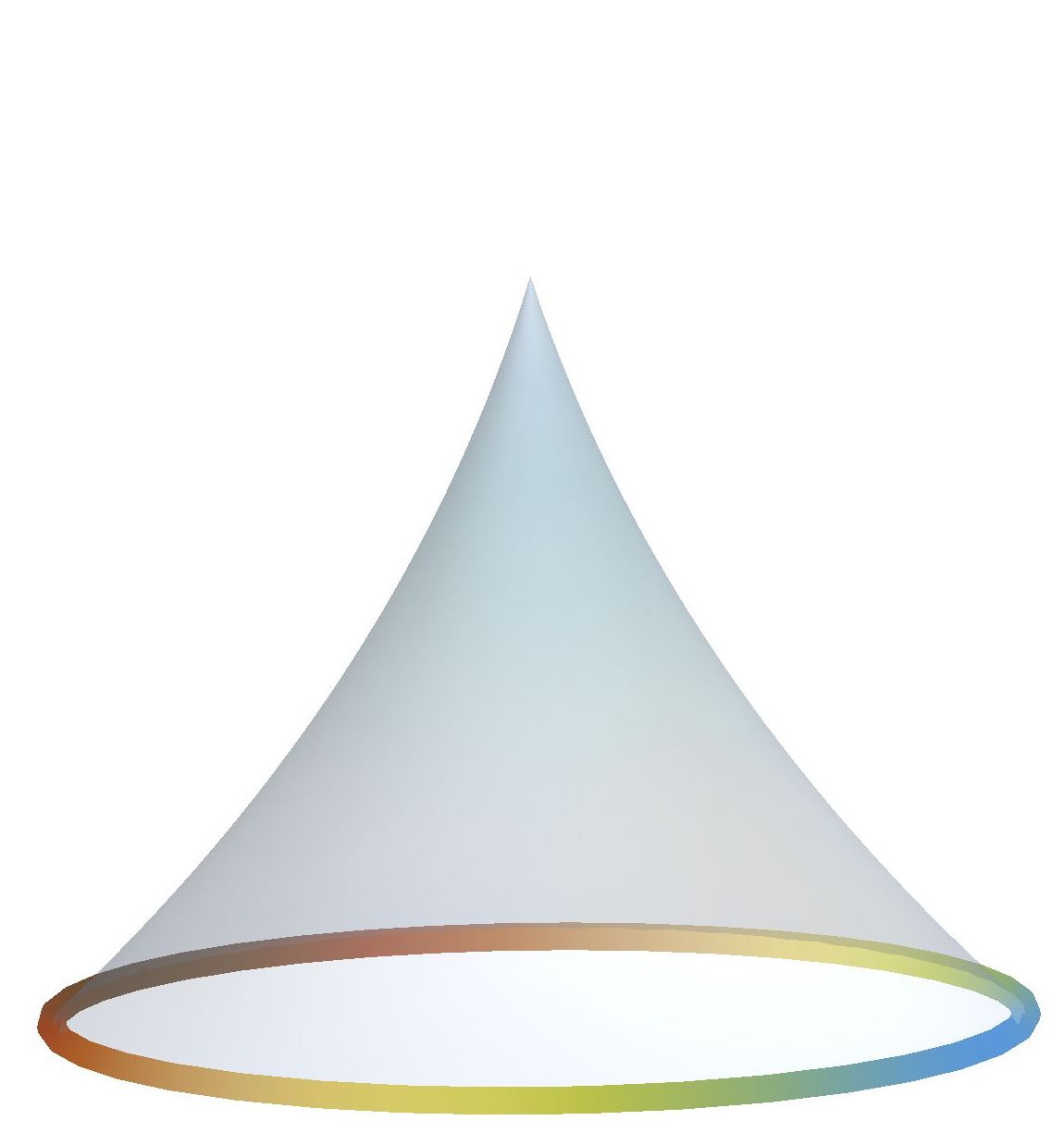}
\end{subfigure}
\begin{subfigure}[b]{0.35\linewidth}
\includegraphics[width=\linewidth]{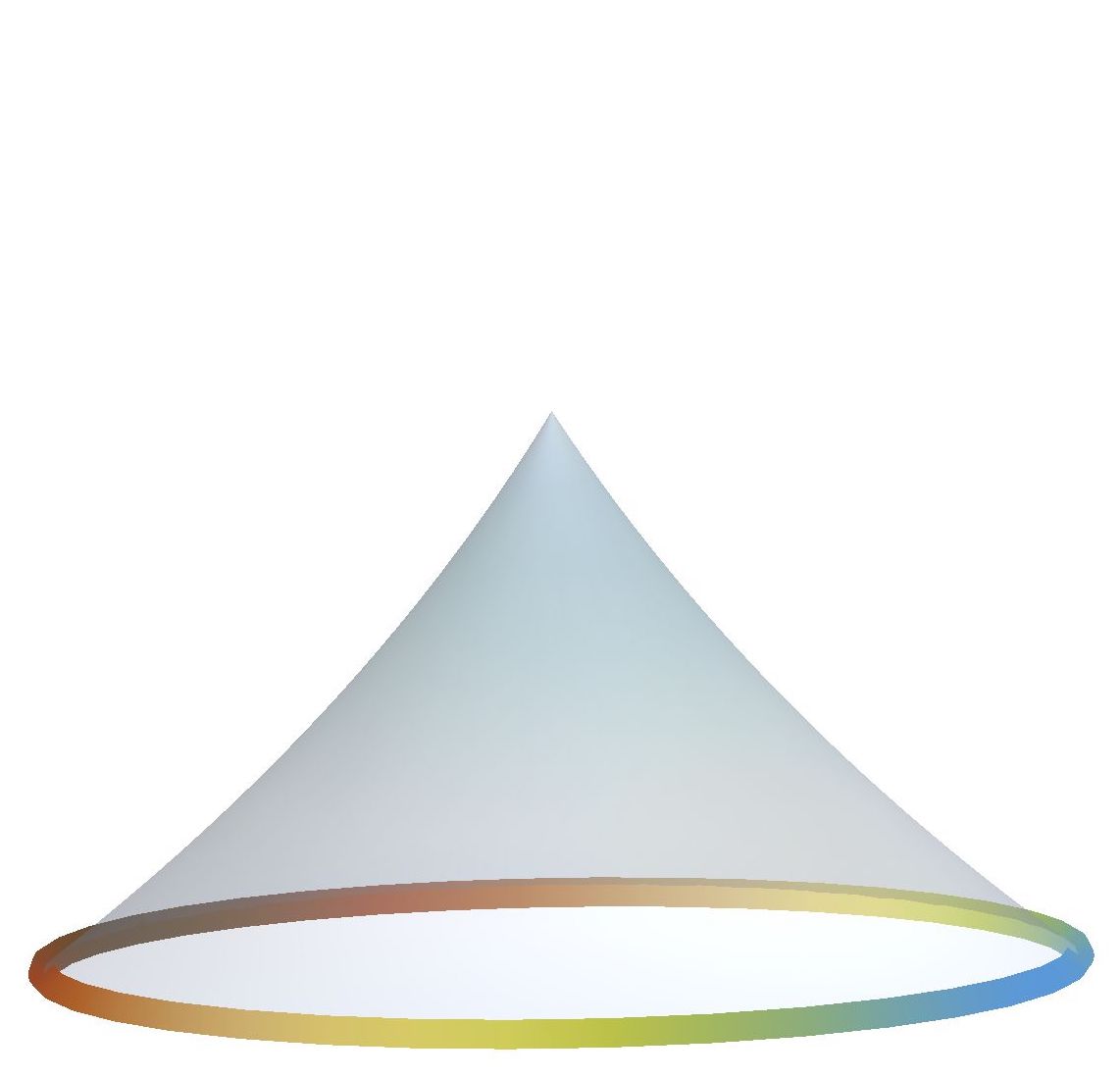}
\end{subfigure}}
\caption{Axially symmetric domains satisfying $K\equiv 2c_o H$, for $c_o\neq 0$. Along the boundary circle $H=2c_o$ holds.}
\label{tents}
\end{figure}

By using \eqref{fichi} to make a change of variable in the integral of above parameterization and considering $\chi\in\left[-c_o,-2c_o/3\right]$, we show, in Figure \ref{tents}, three axially symmetric domains satisfying the Euler-Lagrange equation \eqref{EL11} for $c_o\neq 0$. Moreover, the boundary of these domains is a circle along which $\kappa_n=2c_o$. Therefore, for suitable energy parameters $\sigma$, $\eta$ and $\beta$ such that $\sigma=-4\eta c_o$ and $\kappa\dot{\Lambda}-\Lambda-\beta=0$ holds along the boundary circle, these domains satisfy the Euler-Lagrange equations \eqref{EL11}-\eqref{EL14}.

\begin{rem} These axially symmetric surfaces are only continuous at $r=0$. Due to this lack of regularity, these domains are not solutions of the equation \eqref{fias}.
\end{rem}

\subsection{Constant Mean Curvature Equilibria}

This section will focus on equilibrium configurations for $F[X]$, \eqref{F}, with constant mean curvature. A first approach establishes the following result, which distinguishes between two different classes of constant mean curvature critical domains.

\begin{prop}\label{twoCMC} Let $X:\Sigma\rightarrow\mathbb{R}^3$ be a constant mean curvature $H\equiv H_o$ immersion of disk type which is critical for $F[X]$, \eqref{F}. Then, one of the following holds:
\begin{enumerate}[(i)]
\item If $H\equiv H_o\neq c_o$, the surface is isoparametric and satisfies the relation
$$2(p-1)H^2+2c_o H-pK=0\,,$$
between the mean and Gaussian curvatures, $H$ and $K$, respectively. In this case, either $p=0$ or $\eta\neq 0$.
\item The surface has constant mean curvature $H\equiv H_o=c_o$. This case appears provided that $p\geq 2$ holds.
\end{enumerate}
\end{prop}
\textit{Proof.} Let $X:\Sigma\rightarrow\mathbb{R}^3$ be an immersion with constant mean curvature $H\equiv H_o$. If it is critical for $F[X]$, \eqref{F}, then the Euler-Lagrange equation \eqref{ELF1} must hold on $\Sigma$. For constant mean curvature $H\equiv H_o$, this equation reduces to
$$\left(H_o-c_o\right)^{p-1}\left(2\left[p-1\right]H_o^2-pK+c_oH_o\right)=0\,.$$
Therefore, either $H_o=c_o$ and (ii) holds (note that this case is only possible if $p>1$), or
$$2\left(p-1\right)H_o^2-pK+c_oH_o=0.$$

In the latter, either $p=0$ and Theorem \ref{t} yields that the surface is a planar disk, or the above equation shows that the Gaussian curvature of the surface is also constant, which implies the surface is isoparametric. Moreover, \eqref{ELF2} holds along the boundary, so $\eta\neq 0$ necessarily since $H_o\neq c_o$. This completes (i). \hfill$\square$
\\

For a more refined analysis, first consider the case where the constant mean curvature satisfies $H\equiv H_o\neq c_o$. As a consequence of Proposition \ref{planar}, $p=0$ must hold necessarily for planar disks. This case was studied in detail in Section 4.1. On the other hand, there are non-planar critical disks with $H\equiv H_o\neq c_o$ for arbitrary natural numbers $p$, as is shown in the following result.

\begin{prop}\label{caps} Let $C$ be a circle whose curvature $\kappa\equiv\kappa_o>0$ satisfies $\kappa\dot{\Lambda}-\Lambda-\beta=0$. Then, for any integer $p>0$, there exists a spherical cap spanning $C$ with $H\equiv H_o\neq c_o$ which is critical for $F[X]$, \eqref{F}, provided that $\kappa_o\geq -H_o$ and:
\begin{enumerate}[(i)]
\item If $p=1$, then $c_o<0$, $\eta>0$, $\sigma=-4\eta c_o$ and $H_o=2c_o<0$.
\item If $p=2$, then $c_o=0$, $\eta<0$, $\sigma=-\eta$ and $H_o<0$ is arbitrary.
\item If $p>2$, then $c_o>0$, $\eta\neq 0$, $\sigma p\,^p(-c_o)^{p-2}+4\eta(p-2)^{p-2}=0$ and 
$$H_o=\frac{-2c_o}{p-2}<0\,.$$
In particular, if $p>2$ is even, $\eta<0$; while, if $p>2$ is odd, $\eta>0$.
\end{enumerate}
\end{prop} 
\textit{Proof.} Fix $C$ a circle whose curvature $\kappa\equiv\kappa_o$ satisfies $\kappa\dot{\Lambda}-\Lambda-\beta=0$. Then, there exists a spherical cap spanning $C$, if and only if, the constant mean curvature of the sphere $H\equiv H_o$ satisfies $\kappa_o\geq -H_o$ (recall that, due to our orientation, constant mean curvature surfaces have non positive mean curvature). This first condition is clear, since otherwise the radius of the circle would be bigger than the radius of the sphere and, hence, the circle cannot lie in the sphere. 

Next, we are going to consider the restrictions for which these spherical caps are critical for $F[X]$, \eqref{F}, with $p>0$ and $H\equiv H_o\neq c_o$. From Proposition \ref{twoCMC}, since $H\equiv H_o\neq c_o$, the following equation must hold: 
\begin{equation}\label{Hproof}
2(p-1)H_o^2+2c_o H_o-pK=\left(\left[p-2\right]H_o+2c_o\right)H_o=0\,,
\end{equation}
where we have used that spheres are totally umbilical, i.e. $K=H^2$ holds.

Consider now the case $p=1$. Then, from equation \eqref{Hproof}, we conclude that $H_o=2c_o<0$ must hold. Using this in the boundary condition \eqref{ELF4}, we get that $\sigma=-4\eta c_o$ and, thus, $\eta>0$ necessarily. Observe that the boundary condition \eqref{ELF3} is an identity, while \eqref{ELF2} simplifies to $\kappa_n=2c_o$, which is also an identity. 

Next, assume that $p=2$ holds. From \eqref{Hproof} we obtain that $c_o=0$ while there is no restriction for $H_o$ (apart from $\kappa_o\geq -H_o$ as mentioned at the beginning). We combine this with the boundary conditions \eqref{ELF2}-\eqref{ELF4} (equation \eqref{ELF3} is an identity), obtaining
\begin{eqnarray*}
\sigma H_o+\eta\kappa_n&=&0\,,\\
\sigma H_o^2+\eta K&=&0\,,
\end{eqnarray*}
which, since $H_o=\kappa_n$ and $K=H_o^2$ for a sphere, are both satisfied if and only if $\sigma=-\eta>0$ (c.f. Theorem 3.1 of \cite{Palmer-Pampano-2}). 

Finally, consider the case where $p>2$. For these values, we obtain from \eqref{Hproof} that $H_o=-2c_o/(p-2)$. Moreover, the Euler-Lagrange equations \eqref{ELF2}-\eqref{ELF4} (as above \eqref{ELF3} is an identity) read
\begin{eqnarray*}
\sigma p\left(H_o-c_o\right)^{p-1}+2\eta\kappa_n&=&0\,,\\
\sigma\left(H_o-c_o\right)^p+\eta H_o^2&=&0\,.
\end{eqnarray*}
If we use that $H_o=-2c_o/(p-2)$, it turns out that above equations are the same and that they reduce to $\sigma p\,^p(-c_o)^{p-2}+4\eta(p-2)^{p-2}=0$, as stated. In particular, from this condition we conclude the different possible signs of $\eta$ depending on $p$ being even ($\eta<0$) or odd ($\eta>0$). This finishes the proof. \hfill$\square$
\\

From now on, let us search for constant mean curvature $H\equiv H_o=c_o$ critical domains for the energy $F[X]$, \eqref{F}, when $p\geq 2$. Note that the present choice of orientation (the normal $\nu$ pointing out of a convex domain) implies that the mean curvature $H$ cannot be a positive constant, and therefore $c_o\leq 0$ holds.

\begin{rem} All subsequent results about critical immersions hold for any real number $p\geq 2$.
\end{rem}

Using $H\equiv H_o=c_o$, equations \eqref{ELF1}-\eqref{ELF4} simplify to
\begin{eqnarray}
\eta\kappa_n&=&0\,,\label{bc1}\\
J'\cdot \nu+\eta\tau_g'&=&0\,,\label{bc2}\\
J'\cdot n-\eta\tau_g^2&=&0\,,\label{bc3}
\end{eqnarray}
on the boundary $\partial\Sigma$. Note that the last equation follows from the definition of the Gaussian curvature along $\partial\Sigma$, i.e. $\eta K=\eta\kappa_n\left(2c_o-\kappa_n\right)-\eta\tau_g^2=-\eta\tau_g^2$ holds on $\partial\Sigma$, since $\eta\kappa_n=0$.

It turns out that the boundary must be a solution of above equations, \eqref{bc1}-\eqref{bc3}. This leads to the following characterization.

\begin{thm}\label{crit} Let $X:\Sigma\rightarrow\mathbb{R}^3$ be a constant mean curvature $H\equiv H_o=c_o$ immersion critical for $F[X]$, \eqref{F}, with $p\geq 2$. Then the boundary is a closed critical curve for the energy
\begin{equation}\label{energycurves}
\mathbf{\Theta}_\varpi[C]:=\int_C\left(\Upsilon(\kappa)+\varpi\tau+\beta\right)ds
\end{equation}
where $\Upsilon(\kappa):=\Lambda(\kappa)+\eta\kappa$. Moreover, if the surface is embedded, then the boundary curve must be simple.
\end{thm}
\textit{Proof.} Let the energy $\mathbf{\Theta}_\varpi[C]$, \eqref{energycurves}, act on the space of closed curves. By standard arguments involving integration by parts (see Section 2), a straightforward calculation yields
$$\delta\mathbf{\Theta}_\varpi[C]=\oint_{C}V'\cdot \delta C\,ds$$
where the vector field $V$ is defined along $C$ by
\begin{equation}\label{V}
V:=\left(\kappa\dot{\Lambda}-\Lambda-\beta\right)T+\dot{\Lambda}_sN+\left(\tau\dot{\Lambda}+\eta\tau-\varpi\kappa\right)B\,.
\end{equation}
Hence, $V'= 0$ represents the Euler-Lagrange equation associated to $\mathbf{\Theta}_\varpi[C]$, \eqref{energycurves}.

Now, a constant mean curvature $H\equiv H_o=c_o$ immersion $X:\Sigma\rightarrow\mathbb{R}^3$ critical for $F[X]$, \eqref{F}, with $p\geq 2$ has a boundary curve which satisfies \eqref{bc1}-\eqref{bc3}. A direct check shows that $V'=0$ holds along $\partial\Sigma$, and hence the boundary is a closed critical curve. \hfill$\square$
\\

Critical curves for $\mathbf{\Theta}_\varpi[C]$, \eqref{energycurves}, are the center lines of generalized Kirchhoff elastic rods (c.f. the energy $\mathcal{K}[C]$, \eqref{Kirchhoff}, of Section 2). Observe that the Noether symmetry associated to the translational invariance of $\mathbf{\Theta}_\varpi[C]$ along a critical curve is given by $V\equiv A$, \eqref{V}, for a constant vector $A\in\mathbb{R}^3$. In other words, $V\equiv A$ represents a first integral of the Euler-Lagrange equations associated to $\mathbf{\Theta}_\varpi[C]$, \eqref{energycurves}. Alternatively, this conserved quantity can be obtained directly from \eqref{bc1}-\eqref{bc3}. Indeed, using \eqref{bc1}, equations \eqref{bc2} and \eqref{bc3} can be integrated once, giving
$$V:=J+\eta\tau_g\nu\equiv A\,,$$
for a constant vector $A\in\mathbb{R}^3$. Recall that the vector field $J$ is described in \eqref{J}.

If $A$ is the zero vector, so that $V\equiv 0$ holds along a critical curve $C$, then this critical curve is a circle satisfying $\kappa\dot{\Lambda}-\Lambda-\beta=0$. On the other hand, if $V$ is a nonzero constant vector, then closed critical curves are sought which represent the boundary of a critical disk, therefore some boundary conditions must be satisfied. One of these conditions is a direct consequence of Proposition \ref{rescaling}.

\begin{prop} Let $X:\Sigma\rightarrow\mathbb{R}^3$ be a constant mean curvature $H\equiv H_o=c_o$ immersion critical for $F[X]$, \eqref{F}, with $p\geq 2$. Then, if the boundary is not a circle, it satisfies
$$\int_0^\rho\left(\kappa\dot{\Lambda}-\Lambda-\beta\right)ds=0\,,$$
where $\rho$ represents the period of the curvature of the boundary curve.
\end{prop}
\textit{Proof.} The proof proceeds as in Proposition \ref{intcon}, combining the rescaling condition \eqref{rescalingequation} with $H\equiv H_o=c_o$ on $\Sigma$. \hfill$\square$
\\

This condition can also be deduced from the definition of the vector field $V$. Since the vector field $V$ is constant along the closed boundary $C$, it follows that
$$\oint_C\left(\kappa\dot{\Lambda}-\Lambda-\beta\right)ds=\oint_C J\cdot T\,ds=\oint_C V\cdot T\,ds=\oint_C V\cdot C'\,ds=-\oint_C V'\cdot C\,ds=0\,.$$
Finally, the periodicity of the curvature of $C$ finishes the argument.

Two essentially different cases are now established depending on whether $\eta$ is zero or not.

If $\eta=0$ holds, the boundary conditions are decoupled from the data of the surface interior, and, consequently, there is no restriction on how the surface meets the boundary (see equation \eqref{bc1}). Hence, the initial problem of searching for critical immersions for $F[X]$, \eqref{F}, with $p\geq 2$ has been reduced to seeking a constant mean curvature $H\equiv H_o=c_o$ surface spanning a fixed boundary which is the center line of a generalized Kirchhoff elastic rod, i.e. a  curve critical for $\mathbf{\Theta}_\varpi[C]$, \eqref{energycurves}. That is, the desired surface should solve the volume constrained Plateau Problem, whose solutions for disk type surfaces was given by Hildebrant, \cite{Hildebrant} (see also \cite{Heinz}). Note that these solutions may have branch points.

In order to produce specific examples, we now fix $\Lambda(\kappa)=\kappa^2$ and numerically solve the Plateau problem for the fixed boundary, applying the mean curvature flow. For this choice of boundary energy, boundary curves represent center lines of classical Kirchhoff elastic rods, \cite{Langer-Singer}. Closed center lines of Kirchhoff elastic rods are included in suitable regular homotopy of closed curves. Each class of homotopy contains exactly one (classical) elastic curve (i.e. $\varpi=0$), one self-intersecting center line and one closed curve of constant curvature. Most of these curves are embedded and lie on rotational tori and, hence, they can be understood as torus knots, \cite{Ivey-Singer}. A surface whose boundary is a knot is a Seifert surface. Since we are seeking disk type surfaces, by comparing its genus with the genus of the torus knot, we conclude that boundary curves must be unknotted (for details, see Theorem 4.2 of \cite{Palmer-Pampano-2}). 

\begin{figure}[H]
\makebox[\textwidth][c]{
\centering
\begin{subfigure}[b]{0.3\linewidth}
\includegraphics[width=\linewidth]{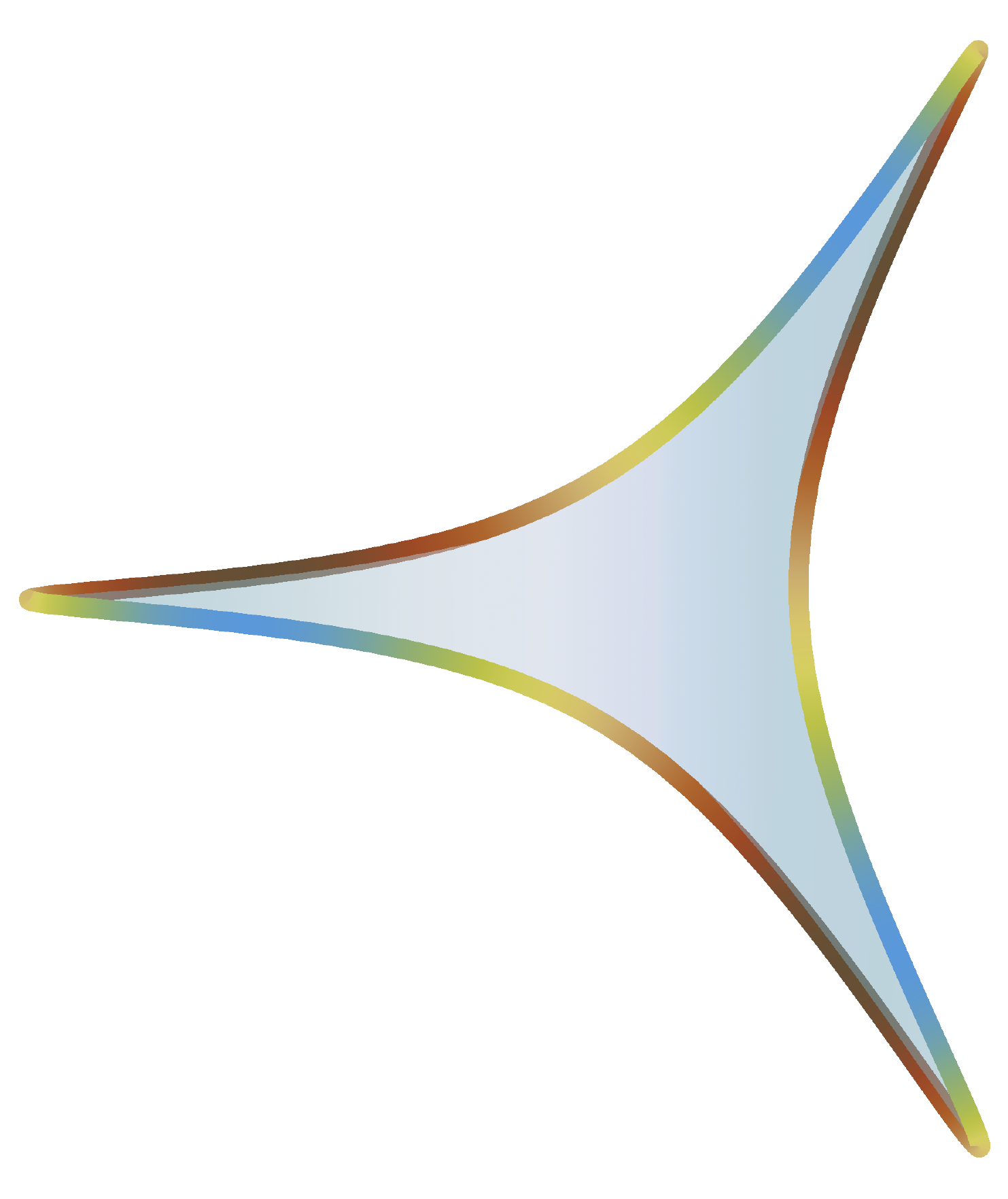}
\end{subfigure}
\begin{subfigure}[b]{0.35\linewidth}
\includegraphics[width=\linewidth]{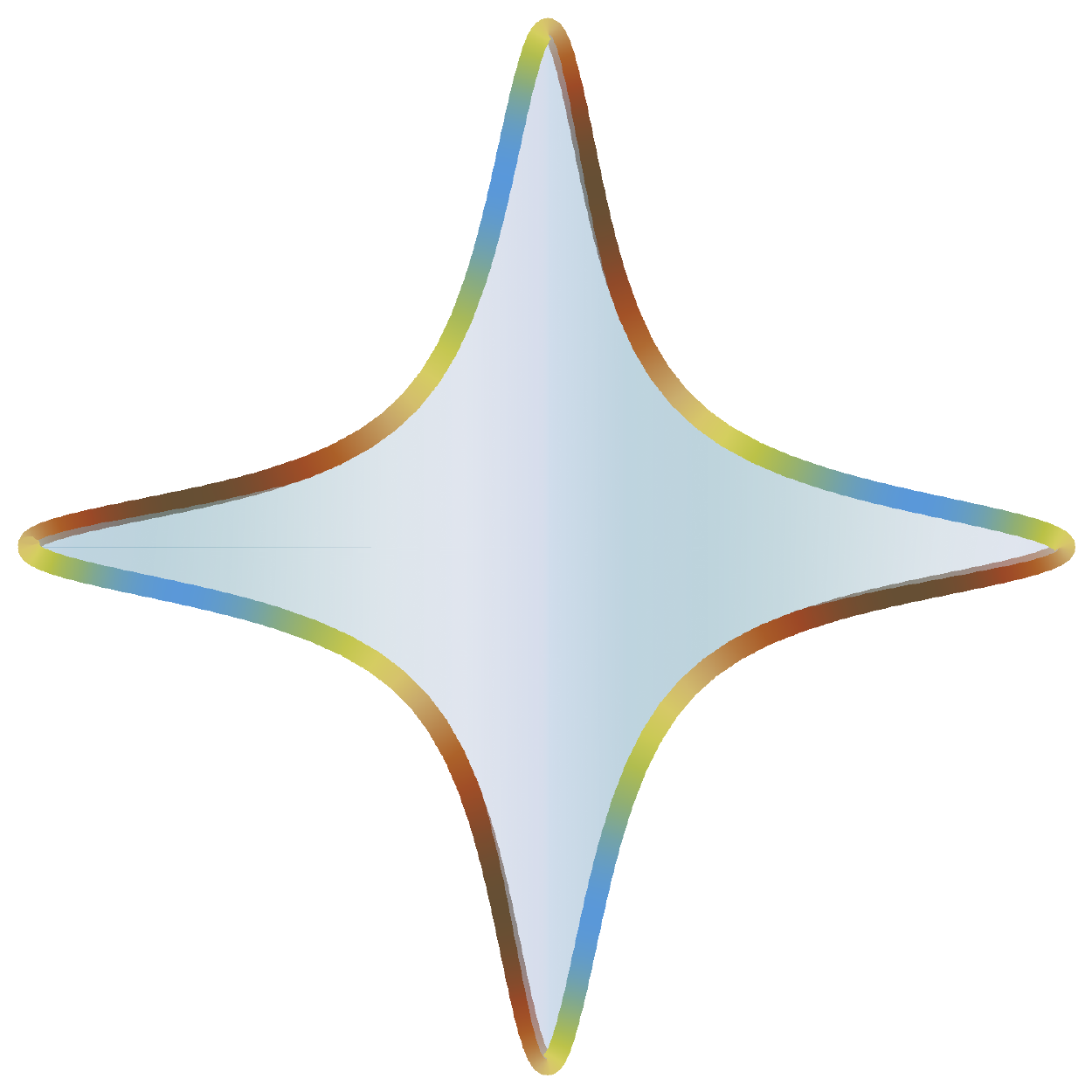}
\end{subfigure}
\,
\begin{subfigure}[b]{0.334\linewidth}
\includegraphics[width=\linewidth]{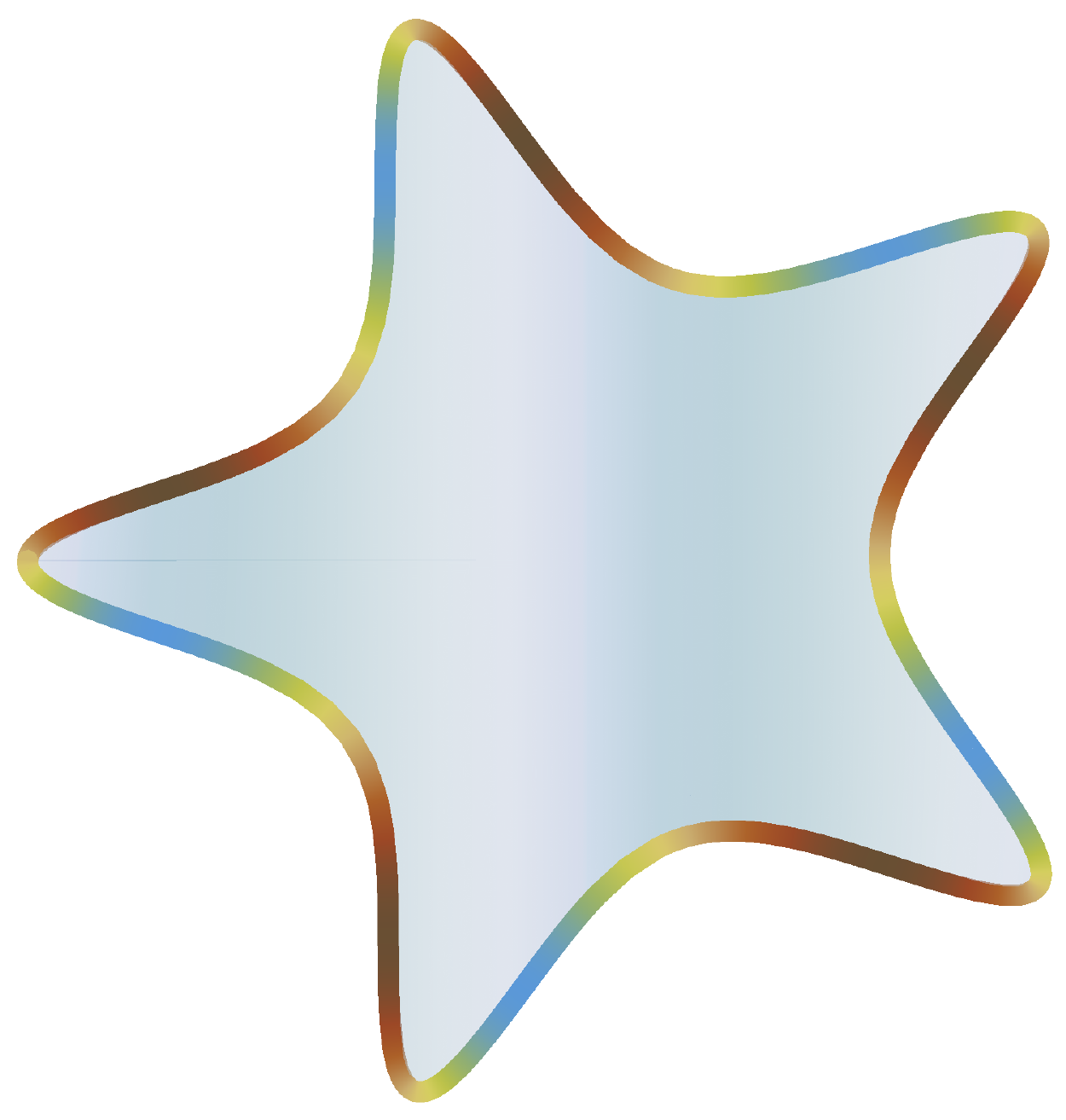}
\end{subfigure}}
\caption{Minimal surfaces of disk type spanned by center lines of Kirchhoff elastic rods representing (unknotted) torus knots. These configurations are critical for the energy $F[X]$, \eqref{F}, for $c_o=0$, $p\geq 2$ and $\eta=0$. The boundary energy for this examples is given by $\Lambda(\kappa)=\kappa^2$, $\varpi\neq 0$ and $\beta>0$.}
\label{CMC}
\end{figure}

Illustrations of minimal surfaces of disk type bounded by elastic curves representing (unknotted) torus knots can be found in \cite{Palmer-Pampano-2}. In Figure \ref{CMC}, we account for twisting, i.e. we consider $\varpi\neq 0$, and show some minimal disk type surfaces bounded by center lines of Kirchhoff elastic rods representing (unknotted) torus knots.

On the other hand, if $\eta\neq 0$, then \eqref{bc1} shows that $\kappa_n\equiv 0$ holds along the boundary, and hence the contact angle between the surface and the center line representing the boundary must satisfy $\theta\equiv \pm\pi/2$. Local examples of such surfaces can be produced using Bj\"{o}rling's formula (see \cite{Bjorling} and \cite{Schwarz}) for $H\equiv H_o=c_o=0$ or the loop group formulation (\cite{Brander-Dorfmeister}) for $H\equiv H_o=c_o\neq 0$. However, in this case, global solutions are always planar disks and $c_o=0$ holds.

\begin{thm}\label{48} Let $X:\Sigma\rightarrow\mathbb{R}^3$ be a constant mean curvature $H\equiv H_o=c_o$ immersion of a topological disk critical for $F[X]$, \eqref{F}, with $p\geq 2$ and $\eta\neq 0$. Then, the surface is planar and $c_o=0$.
\end{thm}
\textit{Proof.} Since $X:\Sigma\rightarrow\mathbb{R}^3$ is a critical immersion for $F[X]$, \eqref{F}, with $p\geq 2$ and $\eta\neq 0$, it follows from \eqref{ELF2} that $\kappa_n\equiv 0$ holds along the boundary.

Next, as done in Theorem \ref{t}, we adapt an argument due to Nitsche (see \cite{Nitsche}) involving the Hopf differential (for details, see again \cite{Gruber-Pampano-Toda} or \cite{Palmer-Pampano-2}) to prove that a disk type constant mean curvature $H\equiv H_o=c_o$ surface satisfying $\kappa_n\equiv 0$ along the boundary is a planar domain. Therefore, $H\equiv H_o=c_o=0$ holds. \hfill$\square$
\\

As proved in Proposition \ref{planar}, the boundary of a planar disk critical for $F[X]$, \eqref{F}, with $p\geq 2$ is a critical curve for
$$\mathbf{\Theta}[C]=\int_C\left(\Lambda(\kappa)+\beta\right)ds\,,$$
satisfying $\varpi\kappa'=0$. Moreover, if the surface is embedded, then the boundary curve is also simple. Observe that this is compatible with the result of Theorem \ref{crit}. In fact, for planar curves $\tau\equiv 0$ and the total curvature is constant on a regular homotopy class.

\begin{rem} The result of Theorem \ref{48} relies heavily on the assumption that the surface is a topological disk. Indeed, there exist constant mean curvature $H\equiv H_o=c_o$ immersions critical for $F[X]$, \eqref{F}, with $p\geq 2$ and satisfying $\kappa_n\equiv 0$ along the boundary for different topological types. For suitable radii of the boundary circles, the nodoidal domains of \cite{Palmer-Pampano-2} are embedded examples of this fact. 
\end{rem}

\section*{Acknowledgements}

We would like to thank Professor E. Aulisa for his help with the numerical computations to obtain the figures of this paper.

\begin{flushleft}
Anthony G{\footnotesize RUBER}\\
Department of Mathematics and Statistics, Texas Tech University, Lubbock, TX, 79409, USA\\
E-mail: anthony.gruber@ttu.edu
\end{flushleft}

\begin{flushleft}
\'Alvaro P{\footnotesize \'AMPANO}\\
Department of Mathematics and Statistics, Texas Tech University, Lubbock, TX, 79409, USA\\
E-mail: alvaro.pampano@ttu.edu
\end{flushleft}

\begin{flushleft}
Magdalena T{\footnotesize ODA}\\
Department of Mathematics and Statistics, Texas Tech University, Lubbock, TX, 79409, USA\\
E-mail: magda.toda@ttu.edu
\end{flushleft}

\end{document}